\documentclass{article}
\usepackage{fancyhdr}
\usepackage{amsmath,amssymb}
\usepackage{array}
\usepackage{mdwmath}
\usepackage{mdwtab}
\usepackage[mathscr]{eucal}
\usepackage{graphicx}
\usepackage{times,pifont}
\usepackage{srcltx}
\usepackage{subfigure}
\usepackage{color}
\usepackage{epsfig}
\usepackage{enumerate}
\usepackage{morefloats}

\newcommand{\ve}{\varepsilon}

\newcommand{\br}{\mathbb{R}}
\newcommand{\bldr}{\boldsymbol{r}}

\newcommand{\url}[1]{\texttt{#1}}

\newcommand{\dis}[2]{\text{\bf dist}_{#2}[#1]}

\newcommand{\spr}[2]{\left\langle #1; #2 \right\rangle}

\newcommand{\blr}{\boldsymbol{r}}

\newtheorem{lem}{Lemma}[section]
\newtheorem{thm}{Theorem}[section]
\newtheorem{remark}{Remark}
\newtheorem{definition}{Definition}
\newtheorem{assumption}{Assumption}
\newtheorem{prop}{Proposition}[section]
\newtheorem{corollary}{Corollary}
\newtheorem{claim}{Claim}
\newtheorem{observation}{Observation}

\newcommand{\pf}{{\it Proof:\;}}
\newcommand{\epf}{$\quad \bullet$}

\begin{document}

\title{Reactive Kinematic Navigation among Moving and Deforming Obstacles with Global Proofs}

\author{Alexey S. Matveev\\almat1712@yahoo.com\\Department of Mathematics and Mechanics,\\Saint Petersburg
University,\\Universitetskii 28, Petrodvoretz,\\St. Petersburg, 198504,
Russia\\
\\
Michael C. Hoy\\mch.hoy@gmail.com\\School of Electrical Engineering and
Telecommunications,\\The University of New South Wales,\\Sydney, 2052,
Australia\\
\\
Andrey V. Savkin\\a.savkin@unsw.edu.au\\School of Electrical Engineering and
Telecommunications,\\The University of New South Wales,\\Sydney, 2052,
Australia\\}

\maketitle









\begin{abstract}                
We present a method for guidance of a Dubins-like vehicle with saturated control towards a target in a steady simply connected maze-like environment. The vehicle always has access to to the target
relative bearing angle (even if the target is behind the obstacle or is far from the vehicle) and the distance to the nearest point of the maze if it is within the given sensor range. The proposed control law is composed by biologically inspired reflex-level rules. Mathematically rigorous analysis of this law is provided; its convergence and performance are
confirmed by computer simulations and experiments with real robots.
\end{abstract}

\section{Introduction}
The capability to successfully operate in dynamic and a priori
unknown environments is an ubiquitous key requirement to mobile
robots. Various scenarios gave rise to a plenty of relevant
algorithms under various sets of assumptions, with the focus being
on collision avoidance. However despite extensive research, this
issue still represents a real challenge in many cases, often because
of numerous uncertainties inherent in the scenario and deficiencies
in a priory knowledge.
\par
With focus on the planning horizon, available rich variety of algorithms can be classified
into global and local planners \cite{LZL07}.
Global planners (GP) generate a complete trajectory based on a
comprehensive model of the scene, which is built from a priori and
sensory data \cite{Latom91}. For dynamic scenes, this approach is
exemplified by several techniques (surveyed in e.g.,
\cite{LaLaSh05,KuVu06}), including state-time space
\cite{ErLoPe87,Fraich99,ReSh94}, velocity obstacles
\cite{FS98,LaLaSh05}, and nonholonomic planners \cite{QuWaPl40}.
Many GP's can be accompanied with firm guarantees of not only
collision avoidance but also achieving a global objective. However
GP's are computationally expensive and hardly suit real-time
implementation; NP-hardness, the mathematical seal for
intractability, was established for even the simplest problems of
dynamic motion planning \cite{Canny88}. A partial cure was offered
in the form of randomized path planning architectures
\cite{HsKiLaRo02,FraDahFe02}. At the same time, all GP's are hardly
troubled by data incompleteness and unpredictability of the dynamic
scene, up to failure in generation of an entire plan.
\par
Local planners (LP) iteratively re-compute only a short-horizon path
on the basis of sensory data about a nearest fraction of the
environment. This reduces the calculation time and dependence on a
priory knowledge. Many of the related techniques, such as the
dynamic window \cite{SeMaPe05,FBTh97}, the curvature velocity
\cite{Simm96}, and the lane curvature \cite{NaSi98} approaches in
fact treat the obstacles as static. On the other side, approaches
like velocity obstacles \cite{FS98}, collision cones \cite{ChGh98},
or inevitable collision states \cite{FraAs03,OwMo06} assume known
and predictable obstacle velocity. However real-world scenarios may
involve obstacles whose predictability ranges from full to none
\cite{LaVa06}. An intermediate approach assumes predictability with
uncertainty and involves a model-based estimation of the region that
may be occupied by an obstacle in the future and that typically
expands as time progresses \cite{LaLaSh05,WuHo12}.
\par
In hardly predictable environments, guarantees provided for LP's up
to now are typically confined to collision avoidance, whereas
achieving the global objective remains an open issue. Moreover with
the exception of \cite{WuHo12}, safety typically concerns only for a
nearest future, whereas its propagation until the very end of the
experiment is not guaranteed.
\par
Previous work has shown navigation schemes based on rigorously proven feedback control laws
have several advantages over alternatives \cite{TeSav10rb,TS10,savkin2013reactive,savkin2013simple,hoy2012collision,MaHoSa13a,MTS11}. This work adopts this design process, however it is applied
to a new problem specification.
\par
The body of the paper is organized as follows. Section~2
describes the system, problem setup, and the proposed navigation strategy.
Assumptions and controller tuning are discussed in Section~3. Section~5 presents the main results. Section~7 is devoted to simulations and experiments with real robots, whereas Section~8 offers brief conclusions. The proofs are given in Appendices~A, B.
\par
\newlength{\mymini}
\setlength{\mymini}{\columnwidth} \addtolength{\mymini}{-5.0pt}
The following notations are adopted in the paper.
\\
$\partial O(t)$--- the boundary of the moving obstacle $O(t)$; \\
$\spr{\cdot}{\cdot}$ --- the standard Euclidian inner product in
$\br^2$;
\\$| \cdot |$ --- the standard Euclidian norm in $\br^3$.
\section{Problem Setup}
\label{S2} We consider a planar point-wise robot traveling in a
two-dimensional workspace. The robot is controlled by the
time-varying linear velocity $\vec{v}$ whose magnitude does not
exceed a given constant $v >0$. The robot's position $\bldr=(x,y)$
is given by its abscissa $x$ and ordinate $y$ in the world frame.
The objective is to move towards the azimuth given by the unit
vector $\vec{f}$. However pure following the azimuth is impossible
since the scene is cluttered with untraversable obstacles $O_1(t),
\ldots, O_N(t), N \geq 1$, and the robot should always be in the
obstacle-free part of the plane $\bldr(t) \not\in \bigcup_i O_i(t)$.
\par
\begin{figure}[h]
\centering
\subfigure[]{\scalebox{0.35}{\includegraphics{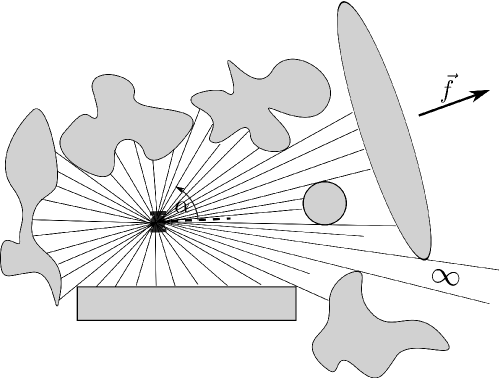}}}
\hfil
\subfigure[]{\scalebox{0.35}{\includegraphics{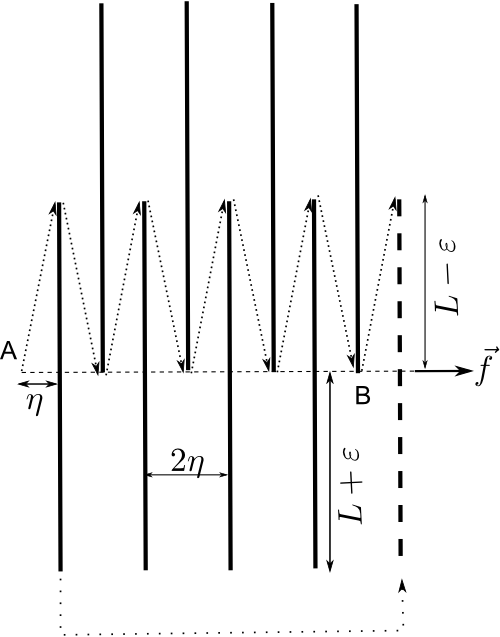}}}
\caption{(a) Robot with a panoramic view; (b) Counterexample} \label{fig.scene}
\end{figure}

The obstacles may undergo general motions, including rotations and
deformations. We require that any obstacle is bounded by a smooth
Jordan curve. Thus we ignore possible inner holes in obstacles on
the ground that they cannot affect the robot's navigation. The
obstacles keep their identities: they do not split into parts, do
not merge together, and do not collide with each other.
\begin{remark}
\rm To take the size of the robot into account, a standard hint is
to treat the obstacle as if it is enlarged and bounded by a proper
equidistant curve. This automatically smoothes the outer cusps of
the obstacle if they do exist. The inner cusps can be smoothed via
approximation to meet the boundary smoothness requirement.
\end{remark}
The robot has a panoramic view on the environment, which view is
shadowed by the obstacles. Specifically we assume a reference frame
attached to the robot. For any polar angle $\alpha$ in this frame,
the robot has access to the distance $d(\alpha,t)$ to the nearest
obstacle in the direction given by $\alpha$ (see
Fig.~\ref{fig.scene}(a)): the point with the relative polar
coordinates $[\alpha,d]$ belongs and does not belong to an obstacle
if $d=d(\alpha,t)$ and $0 \leq d < d(\alpha,t)$, respectively.
If there is no obstacle in the examined direction,
$d(\alpha,t):=\infty$. The desired direction of motion $\vec{f}$ is
also accessed via its relative polar angle $\beta(t)$.
\par
The robot is able to classify enroute obstacles into several
disjoint classes $C_1, \ldots, C_K$, for example, into
\textquoteleft definitely static\textquoteright\, and \textquoteleft
probably moving\textquoteright or into \textquoteleft low
speed\textquoteright\, and \textquoteleft high
speed\textquoteright\, obstacles, etc. This setup covers both the
case of no classification capacity $K=1$, where the robot is
enforced to treat all obstacles in a common way, and the opposite
case where every class contains only one obstacle, which means
capability to recognize the obstacles.
\par
By using these data, the robot should move along the given azimuth.
Unlike many papers in the area, we do not assume that the obstacles
obey a known law of motion. They may move arbitrarily with bounded
speeds. The only assumption is that the robot is faster than the
obstacles, whose necessity will be justified by Lemma~\ref{nec.lem}.
\par
We do not provide specifications of the robot-fixed frame or its
evolution over time since these details does not matter in both
proposed control law and its analysis.

\section{The navigation algorithm}
In what follows, we treat angle $\gamma$ as a cyclic $\gamma \pm 2
\pi = \gamma$ variable associated with the point $\vec{e}(\gamma)$
on the unit circle.
\par
At time $t$, the scene is given by the function $d(\alpha) =
d(\alpha,t)$ of the polar angle $\alpha$. This function and the
angle $\beta = \beta(t)$ of $\vec{f}$ are converted into control
$\vec{v}(t)$ via the following steps.
\begin{itemize}
\item {\it Computation of facets.}
The discontinuity points $\alpha_j$ of $d(\alpha)$ are
determined.\footnote{Since the obstacles have no inner holes
where the robot might appear, such points do exist.} They
partition the circle into maximal arcs $A_k = (\alpha_k^-,
\alpha^+_k)$ on which $d(\cdot)$ is continuous.
If $d(\alpha)=\infty$ somewhere on such arc, it is equal to
$\infty$ on the entire arc; these arcs are dismissed. Every of
the remaining arcs can be assigned to a visible facet of an
obstacle.
\item {\it Enlargement of the facets.}
To bypass moving obstacles, the angular image of any facet is
enlarged so that the closer the facet, the more it is extended.
To this end, a continuous decaying function of the distance $d
\geq 0$
\begin{equation}
\label{eq.delta}
\Delta_i(d) \in \left[ 0,\frac{\pi}{2} \right), \quad \Delta_i(d_1) \geq \Delta_i(d_2)\;
\forall d_1 \leq d_2
\end{equation}
is pre-specified for any class $C_i$ of obstacles. For any facet
$k$, the distance $d_k^{\,\min}:= \inf_{\alpha \in A_k}
d(\alpha)$ to it and the class $i$ of the concerned obstacle are
determined. The {\it extended facet} is given by its range
$\widehat{A}_k $ and {\it extended profile} $d_k(\cdot)$
\begin{gather}
\label{eq.exfac} \widehat{A}_k := \left[ \alpha_k^- - \delta_k,
\alpha^+_k + \delta_k\right], \qquad
\delta_k:=\Delta_i(d_k^{\,\min}),
\\
\label{eq.exfac1} d_k(\alpha) := d [p^{A_k}(\alpha)],
\end{gather}
where $p^{A_k}(\alpha)$ is the point of $A_k$ nearest to
$\alpha$.\footnote{The ranges $\widehat{A}_k$ may overlap.}
\item{\it Generation of the control.} There may be two cases.
\begin{itemize}
\item {\it If desired direction $\beta$ is not obstructed by extended
facets} $\beta \not\in \bigcup_k \widehat{A}_k$, the robot
is driven in this direction
$$
\vec{v} := v \vec{e}(\beta).
$$
\item {\it The desired direction is obstructed}.
An obstructing facet $k$ with the minimal value of
$d_k(\beta)$ is chosen. The end-points $\alpha \in
\widehat{A}_k$ of all extended facets $j$ for which
$d_j(\alpha) \leq d_k(\alpha)$ are gathered in the set
$E_k$. Among the points $\alpha_{\circlearrowleft}$ and
$\alpha_{\circlearrowright}$ of $E_k$ that are
counter-clockwise and clockwise closest to $\beta$,
respectively, we pick the minimizer $\alpha_0$ of the
discrepancy $|\beta -
\alpha_0|$\footnote{$\alpha_{\circlearrowleft}$ in the case
of equal discrepancies} and put
\begin{equation}
\label{c.rule2}
\vec{v} := v \vec{e}(\alpha_0).
\end{equation}
\end{itemize}
\end{itemize}
\par
Though $\Delta_i(\cdot)$ can be chosen constant, decaying functions
allow the robot not to take overly precautionary measures against
faraway obstacles. A reasonable modification results from putting
\eqref{c.rule2} in use only under the extra condition: the distance
to the nearest obstacle is less than a pre-specified threshold
$d_\ast$. This modification will be addressed in ???
\par
The proposed law does not estimate the velocity of the obstacle and
does not attempt to predict its future positions. This is partly
motivated by relevant troubles, up to infeasibility, possible low
accuracy and essential extra computational burden. These problems
are drastically enhanced if obstacles undergo general motions,
including deformations, since then the velocity and future position
are often characterized by infinitely many parameters. Another
reason is that the objective of this research is to put a landmark
that can be reached without the aid of estimation or prediction.
\par
Since the proposed navigation law is discontinuous, the closed-loop
solutions are meant in the Filippov's sense \cite{FIL88}. Given the
initial state, such solution exists and does not blow up in a finite
time since the controls are bounded. Its uniqueness is a more
delicate problem. So we shall address all \textquoteleft
Filippov's\textquoteright solutions with only one exception. In the
situation of repulsive discontinuity surface (the vector field
points away from the surface on both sides), solutions that slide
over the surface are dismissed. The theoretical reason is that they
are unstable and so inviable in the face of unavoidable small
disturbances. The practical reason is that we tacitly assume
sampled-data control by a digital device, which picks at random one
of two control options on the above discontinuity surface and does
not alter it during the sampling time $\tau >0$. This causes
immediate escape from the surface and gives rise to no sliding
solutions as $\tau \to 0$, which is the limit case that is caught by
the model at hand.

\section{Collision avoidance}
The obstacles may not only arbitrarily move but also twist, skew,
wriggle, or be otherwise deformed. We assume that all these are
performed smoothly. To formally describe this, we use the Lagrangian
formalism \cite{Spencer04} by introducing a {\it reference
configuration} $O_{j,\ast} \subset \br^2$ of the $j$th obstacle and
the {\it configuration map} $\Phi_j(\cdot,t): \br^2 \to \br^2$ that
transforms $O_{j,\ast}$ into the current configuration $O_j(t) =
\Phi[O_{j,\ast},t]$.
\begin{assumption} \label{ass}
The reference configuration $O_{j,\ast}$ of any obstacle is compact
and bounded by a smooth Jordan curve. The configuration map
$\Phi_j(\cdot,t)$ is defined on an open neighborhood of $O_{j,\ast}$
and is smooth and one-to-one, the determinant of its Jacobian matrix
is everywhere nonzero.
\end{assumption}
Let $\vec{V}_j(r,t)$ stands for the velocity of the point $\blr \in
O_j(t)$, and $[\vec{T}(\blr,t),\vec{N}(\blr,t)]$ be the Frenet frame
of $\partial O_j(t)$ at $\blr$, respectively. (The set $O_j(t)$ is
to the left to the tangent $\vec{T}$, the normal $\vec{N}$ is
directed inwards $D(t)$.) Finally, $W^T$ and $W^N$ are the
tangential and normal components of the vector $\vec{W}$.
\par
We start with conditions necessary for the robot to be capable of
collision avoidance. The following definition tacitly assumes
scenarios where the avoidance maneuver may start in an arbitrarily
tight proximity of the obstacle.
\begin{definition}
\label{def.av} The robot {\em is capable of avoiding collisions with
the obstacle} $O_j(t)$ if for any initial time $t_0$ and state of
the robot $\bldr(t_0) \not\in O_j(t_0)$, there exists an admissible
velocity profile $|\vec{v}(t)| \leq v, t \geq t_0$ under which the
robot does not collide with this obstacle for $t \geq t_0$.
\end{definition}
Let $[a]_-:= \min \{0, -a\}$ stand for the negative part of $a\in
\br$.
\begin{lem}
\label{nec.lem} Suppose that the robot is capable of avoiding
collisions with the obstacle $O_j(t)$. Then at any time $t$, the
negative part $\left[V_{j}^N(\bldr,t)\right]_-$ of the normal
velocity (i.e., that responsible for the outward direction) of any
boundary point $\bldr \in \partial O_j(t)$ does not exceed the
robot's maximal speed $v$:
\begin{equation}
\label{nec.cond}
\left[V_{j}^N(\bldr,t)\right]_- \leq v .
\end{equation}
\end{lem}
The proofs of the claims stated in this section are given in
Appendix~\ref{app.2}. By the following theorem, putting $<$ in place
of $\leq$ makes the necessary condition \eqref{nec.cond} sufficient,
and under this sufficient and \textquoteleft almost
necessary\textquoteright\, condition, collision avoidance is ensured
by the proposed control law.
\begin{thm}
\label{th.main} Suppose that the negative part of the normal speed
$[V_{j}^N(\bldr,t)]_-$ of any boundary point $\bldr \in O_j(t)$ of
any obstacle is less than the robot's maximal speed $v$, and this
relationship does not degenerate as time goes to $\infty$:
\begin{equation}
\label{spee.req}
\left[ V_{j}^N(\bldr,t)\right]_- \leq v_o^i < v \quad \forall \bldr \in \partial O_j(t), j \in C_i, i, t.
\end{equation}
Then the proposed control law keeps the robot in the obstacle-free
part of the scene if $\Delta_i(\cdot)$ are chosen so that
\begin{equation}
\label{safety.req}
\Delta_i(0)>\arcsin v_o^i/v .
\end{equation}
\end{thm}
Due to \eqref{spee.req}, this is possible. If the speed bound
$v_o^i$ is known, \eqref{safety.req} can be directly used for
controller tuning. Otherwise it provides general guidelines for
tuning: by picking $\Delta_i(0)< \pi/2$ close enough to $\pi/2$,
safety is guaranteed.
\section{Achieving of the main control objective}
In the general setting with non-convex obstacles, this objective
encompasses solving of dynamic mazes. Even for static mazes, this
constitutes a challenging and separate research topic, whereas
deforming mazes lie in the uncharted territory. To reduce the
challenge to a reasonable level, we restrict ourselves to the case
where the obstacles are treated (maybe via approximation) as convex
compact bodies. They still may undergo translations, rotations, and
deformations. We show that the robot with \textquoteleft collision
avoidance\textquoteright capacity \eqref{spee.req},
\eqref{safety.req} achieves the control objective if spacing between
the obstacles or the robot's speed is large enough. Specifically,
the necessary spacing depends on the sizes of the obstacles and the
speed excess ratios $v/v_0^i$ so that the smaller the size or the
larger the ratio, the smaller the spacing. In particular, the
spacing requirement annihilates as $v/v_0^i \to \infty\; \forall$.
\subsection{Counterexample}
\label{subsect.counter} To start with, we show that spacing and size
are really important. Let the robot face $2N$ straight line segments
of length $2L$ perpendicular to $\vec{f}$; see
Fig.~\ref{fig.scene}(b), where $N=4$. They are separated by the
space $2\eta$ and slowly drift in the direction of $-\vec{f}$ with
the speed $v_o<v$. The odd segments are vertically aligned, the even
segments are shifted upwards by $L+\ve$, where $\ve>0, \ve \approx
0$. Initially the robot is in the location $A$. The upper dotted
line illustrates the path of the robot driven by the proposed
control law.
\par
To reach point B, the robot spends no less than $2N(L-\ve)/v$ time
units. For this time, B moves a distance $\geq 2N(L-\ve)v_o/v$ and
reaches a location to the left of $A$ if $2N(L-\ve)v_o/v
>2N\eta \Leftrightarrow e:= v/v_o < (L-\ve)/\eta$. So for any speed excess ratio $e$,
there exist so tight spacing $\eta$ and large obstacle size $L$ that
the overall drift of the robot is contrary to the desired direction.
What is worse, that time is enough for the first segment to arrive
at the dashed position over the lower dotted path if $(2L+2\eta
N)/v_o< 2N(L-\ve)/v \Leftrightarrow e< \left(\frac{L}{N(L-\ve)} +
\frac{\eta}{L-\ve} \right)^{-1}$, which can ensured for any $e$ by
picking $\eta \approx 0$ and $N \approx \infty$. So if every segment
performs such rearrangement just after the robot bypasses it, the
overall drift of the robot contrary to $\vec{f}$ will be never
terminated.
\par
Thus apart from the speed ratio, criteria for achieving of the
control objective should somehow impose requirements on the
obstacles sizes and spacing between them.
\subsection{Main Results}
To flesh out this, we introduce geometric objects, called {\it
hats}. Being the parts of the plane that should be free from the
obstacles, they characterize the required inter-obstacle spacing.
These objects depend on an angle $\delta \in (0,\pi/2)$, which will
depend on the speed ratio in the results to follow.
\par
\begin{figure}[h]
\centering
\subfigure[]{\scalebox{0.25}{\includegraphics{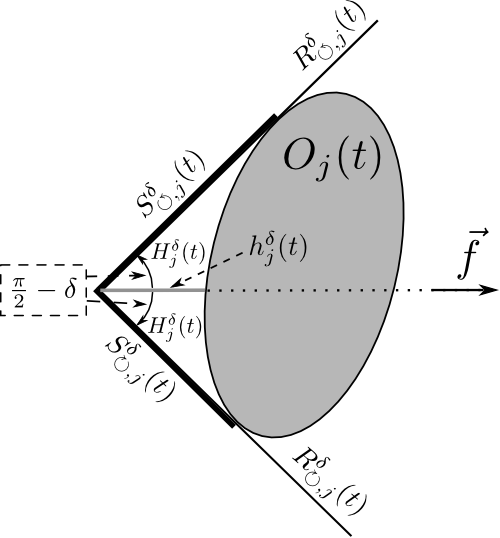}}}
\hfil
\subfigure[]{\scalebox{0.25}{\includegraphics{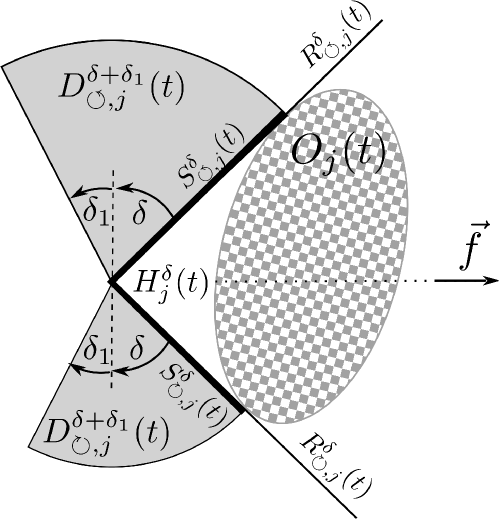}}}
\caption{(a) The $\delta$-hat $H_j^\delta(t)$; (b) The extended hat.}
\label{fig.hat}
\end{figure}
Consider the angle made of two rays that are obtained from the ray
spanned by $\vec{f}$ via the counter clockwise and clockwise
rotation, respectively, though the angle $\pi/2 - \delta$. Let us
translate this angle into the unique position
$\sphericalangle_j^\delta(t)$ in the plane where its both rays
$R_{\circlearrowleft,j}^\delta(t)$ and
$R_{\circlearrowright,j}^\delta(t)$ are tangential to the boundary
of $O_j(t)$; see Fig.~\ref{fig.hat}(a). The {\it $\delta$-hat}
$H_j^\delta(t)$ of the $j$th obstacle at time $t$ is the part of the
plane bounded by
\begin{enumerate}[{\bf 1)}]
\label{semg}
\item the segment $S_{\circlearrowleft,j}^\delta(t)$ of $R_{\circlearrowleft,j}^\delta(t)$
between the vertex of $\sphericalangle_j^\delta(t)$ and the
furthest point of $R_{\circlearrowleft,j}^\delta(t) \cap
O_j(t)$;
\item the segment $S_{\circlearrowright,j}^\delta(t)$ of $R_{\circlearrowright,j}^\delta(t)$ between the
vertex of $\sphericalangle_j^\delta(t)$ and the furthest point
of $R_{\circlearrowright,j}^\delta(t) \cap O_j(t)$;
\item the least part of the boundary $\partial O_j(t)$ that connects these segments
and is closer to the angle vertex.
\end{enumerate}
The {\it height} $h_j^\delta(t)$ of this hat is the distance from
the vertex of $\sphericalangle_j^\delta(t)$ to the obstacle. The
{\it $\delta_1$-extended $\delta$-hat} $
\widehat{H}_j^{\delta,\delta_1}(t)$ is the union of $H_j^\delta(t)$
with the following two sets (see Fig.~\ref{fig.hat}(b)):
\begin{itemize}
\item the sector $D_{\circlearrowleft,j}^{\delta+\delta_1}(t)$
swept by the segment $S_{\circlearrowleft,j}^\delta(t)$ when
counter-clockwise rotating through the angle $ \delta+\delta_1$;
\item the sector $D_{\circlearrowright,j}^{\delta+\delta_1}(t)$ swept by the segment $S_{\circlearrowright,j}^\delta(t)$ when
clockwise rotating through the angle $\delta+\delta_1$.
\end{itemize}
\begin{thm}
\label{th.reach} Suppose that the \textquoteleft
two-sided\textquoteright extension of \eqref{spee.req}
\begin{equation}
\label{spee.req1}
\left| V_{j}^N(\bldr,t)\right| \leq v_o^i < v \quad \forall \bldr \in \partial O_j(t), j \in C_i, i, t
\end{equation}
holds, the obstacles are always convex and do not collide with each
other, and for $\delta_i^\star := \arcsin \frac{v_o^i}{v}$,
\begin{enumerate}[{\bf a)}]
\item the $\delta_j^\star$-extended $\delta_i^\star$-hat of any obstacle of class $C_i$ is always disjoint with
the obstacles of class $C_j$ and
\item for $t=0$ and any $i$, the location $\bldr(0)$ of the robot does not lie in the $\delta_i^\star$-hat of any
obstacle of class $C_i$.
\end{enumerate}
Moreover let {\bf a)} and {\bf b)} be true for a larger value of
$\delta_i^\star$:
\begin{equation}
\label{delta.ineq}
\delta_i^\star
\in \left( \overline{\delta}_i, \frac{\pi}{2}\right), \quad \text{\rm where} \; \overline{\delta}_i := \arcsin \frac{v_o^i}{v}.
\end{equation}
Then the proposed controller drives the robot through the
obstacle-free part of the environment. Moreover, it can be tuned so
that the robot always drifts in the right direction
\begin{equation}
\label{drift}
\spr{\vec{v}(t)}{\vec{f}} \geq \eta >0 \qquad \forall t.
\end{equation}
\par
Specifically, this is true if $\Delta_i(\cdot)$ in \eqref{eq.delta}
is chosen so that
\begin{equation}
\label{delta.range}
\Delta_i\left[ h_j^{\delta_i^\star}(t) \right] = \Delta_i[0] \in
\left( \overline{\delta}_i, \delta_i^\star \right) \quad \forall t, i, j \in C_i.
\end{equation}
\end{thm}
\par
The proof of this theorem is given in Appendix~\ref{app.1}.
\begin{remark}
\label{main.rem} \rm {\bf i)} If b) is true for $\delta_i^\star :=
\overline{\delta}_i$, it also holds for $\delta^\star_i$ satisfying
\eqref{delta.ineq} by d) in Observation~\ref{obs1}. So the concerned
assumption of Theorem~\ref{th.reach} is always valid and is stated
only to introduce $\delta^\star_i$ in \eqref{delta.range}.
\par
{\bf ii)} Theorem~\ref{th.reach} assumes undefined end of the
experiment and so extends the conclusion on all the future. It
remains true if time is confined to a finite interval $[t_0,t_1]$ in
the assumptions and conclusion. Then it suffices to check a) only
for $\delta_i^\star := \overline{\delta}_i$ since this implies a)
with a larger $\delta^\star_i$, like in i).
\par
{\bf iii)} Function $\Delta_i(\cdot)$ satisfying \eqref{delta.range}
does exist. For example given an estimate of the hat height
$h_j^{\delta_i^\star}(t) \leq h_i \; \forall j \in C_i, t$, it
suffices to pick $\Delta_i[0] \in \left( \overline{\delta}_i,
\delta_i^\star \right)$, put $\Delta_i(d) := \Delta_i(0)$ for $d \in
[0,h_i]$, and extend $\Delta_i(\cdot)$ on $[0,\infty)$ as a
continuous decaying function. Another option is a constant
$\Delta_i(\cdot)$, which does not need estimation of the hat height.
\par
{\bf iv)} A guideline given by \eqref{delta.range} is to pick
$\Delta_i(0)$ close to $\overline{\delta}_i$, which may be used if
estimation of $\delta_i^\star$ is troublesome.
\par
{\bf v)} As $v_o^i/v \to 0$, all hats degenerate into parts of the
obstacle boundary. Then a) and b) come to disjointness of the
obstacles and $\bldr(0) \not\in O_j(0)\;\forall j$, whereas
\eqref{delta.range} come to $\Delta_i(0)>0$.
\end{remark}

\section{Illustrations of the main results for special scenarios.}
These scenarios involve obstacles of special shapes. For each
scenario, we start with explicit computation of the concerned hats
and related specifications of Theorem~\ref{th.reach}. In some cases,
these hats appear to be rather complex. So we finish with easier
comprehensible corollaries that result from replacement of the hats
with their upper estimates by simple shapes.
\subsection{Disk-shaped obstacles}
\begin{figure}[h]
\centering
\subfigure[]{\scalebox{0.25}{\includegraphics{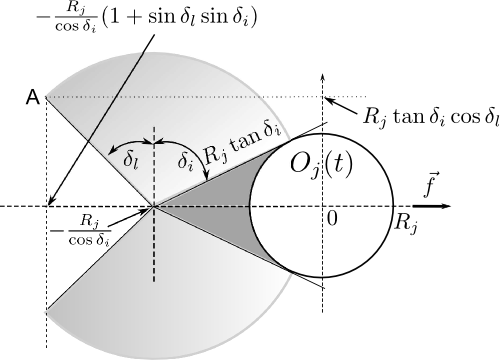}}}
\hfil
\subfigure[]{\scalebox{0.25}{\includegraphics{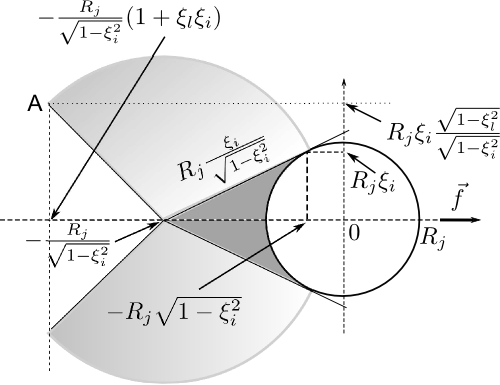}}}
\caption{The hat and extended hat of the disk.}
\label{fig.fac3}
\end{figure}

In robotics research, obstacles are often represented by disks. Let
the $j$th obstacle be the disk of the constant radius $R_j$ with the
center $\bldr_j(t)$ whose velocity is $\vec{v}_j(t)$. The
\textquoteleft collision avoidance\textquoteright\, condition
\eqref{spee.req} is equivalent to \eqref{spee.req1} and shapes into
\begin{equation}
\label{v.d}
|\vec{v}_j(t)| \leq v_o^i <v \qquad \forall j \in C_i.
\end{equation}
Elementary geometrical considerations show that the $\delta_i$-hat
of the $j$th obstacle is the set shadowed in dark in
Fig.~\ref{fig.fac3}(a). The $\delta_l$-extended $\delta_i$-hat is
obtained via union of the hat with two slightly shadowed disk
sectors. In the basic case where $\delta_r=\overline{\delta}_r$,
Fig.~\ref{fig.fac3}(a) shapes into Fig.~\ref{fig.fac3}(b), where
$\xi_r:=v^r_o/v$.
\par
By bringing the pieces together, we arrive at the following.
\begin{prop}
\label{prop.disk} If \eqref{v.d} and \eqref{safety.req} are true,
the robot is driven through the obstacle-free part of the
environment. Suppose also that for any obstacle $j \in C_i$, its
extended hat from Fig.~{\rm \ref{fig.fac3}(b)} is always disjoint
with the obstacles of class $C_l$, and the robot is initially
outside its hat from Fig.~{\rm \ref{fig.fac3}(b)}. Moreover let
these claims be true for the hats from Fig.~{\rm \ref{fig.fac3}(a)}
with some $\delta_r > \overline{\delta}_r \;\forall r$. Then the
robot always drifts in the right direction, i.e., \eqref{drift} is
true, if $\Delta_i(\cdot)$ are chosen so that
\begin{equation}
\label{specif.cont}
\Delta_i\left[ R_j \frac{1- \cos \delta_i}{\cos
\delta_i} \right] = \Delta_i[0] \in
(\overline{\delta}_i,\delta_i).
\end{equation}
\end{prop}
\par
In \eqref{specif.cont}, $R_j$ can be replaced by any its known upper
bound. If this bound is not available, $\Delta_i(\cdot)$ can be
picked constant.
\par
Now we apply Proposition~\ref{prop.disk} to more special scenarios
and, by sacrificing a part of its content, provide convergence
criteria in terms of distances instead of hats.
\par
{\bf A scene cluttered with irregularly and unpredictably moving
disk-shaped obstacles with a common velocity bound} $v_o^i = v_o
<v\; \forall i$. Let $L_{j,k}=L_{k,j}$ be a time-invariant lower
estimate of the spacing between the $j$th and $k$th obstacles. Since
$\delta$-hat is contained by the disk passing through its vertex
(see Lemma~\ref{lem.maxdis} below) and the extended hat is contained
by the disk centered at the origin and passing through point A in
Fig.~{\rm \ref{fig.fac3}(b)}, we arrive at the following claim.
\begin{claim}
The proposed controller guarantees no collisions with the obstacles
and constant drift in the right direction if
\begin{multline}
\label{basic.disk} \frac{\dis{\bldr(0)}{O_j(0)}}{R_j} >
\Omega(\xi):= \left[\frac{1}{\sqrt{1-\xi^2}} -1 \right],
\quad \xi:= \frac{v_o}{v} ; \\
\frac{L_{j,k}}{\max\{R_k,R_j\}} > \Upsilon(\xi):=
\frac{\sqrt{1+3\xi^2}-\sqrt{1-\xi^2}}{\sqrt{1-\xi^2}} ,
\end{multline}
 and the 
\textquoteleft parameter\textquoteright $\Delta(\cdot)$ of the
controller is chosen so that
\begin{equation}
\label{eq.vein}
\Delta\left[ \Omega(\xi) R + \ve \right] = \Delta[0] \in \left( \arcsin \xi, \arcsin \xi^\star \right).
\end{equation}
Here $R \geq R_j\forall j$, $\ve>0$ may be arbitrarily small and
$\xi^\star
> \xi$ is picked so that \eqref{basic.disk} remains true with $\xi$
replaced by $\xi^\star$.
\end{claim}
\par
\begin{figure}[h]
\renewcommand{\arraystretch}{1.2}
\centering
\subfigure[]{\scalebox{0.4}{\includegraphics{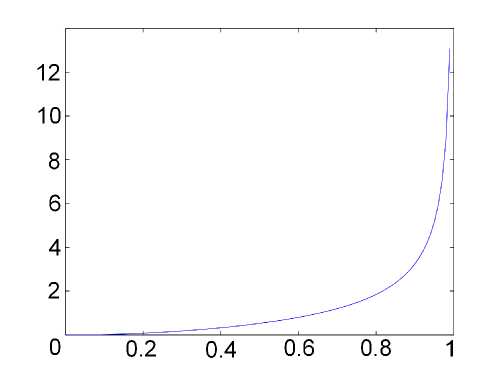}}}
\hfil
\subfigure[]{\raisebox{30.0pt}{\begin{tabular}{|c|c|c|c|c|}
\hline
\tiny $\xi$ & \scriptsize 0 & \tiny 1/8 &\tiny 1/7 &\tiny 1/6
\\
\hline
\tiny $\Upsilon(\xi)$ &\tiny 0 &\tiny 3\% &\tiny 4\%  &\tiny  6\%
\\
\hline
\tiny $\xi$ & \tiny 1/4 & \tiny 1/3 &\tiny 1/2 &\tiny $1/\sqrt{2}$
\\
\hline
\tiny $\Upsilon(\xi)$ &\tiny 12\%     &\tiny 22\% &\tiny 52\%  &\tiny  124\%
\\
\hline
\end{tabular}}}
\caption{(a) Graph of $\Upsilon(\cdot)$; (b) Table of values.}
\label{fig.ups}
\end{figure}
\begin{figure}[h]
\renewcommand{\arraystretch}{1.2}
\centering
\subfigure[]{\scalebox{0.4}{\includegraphics{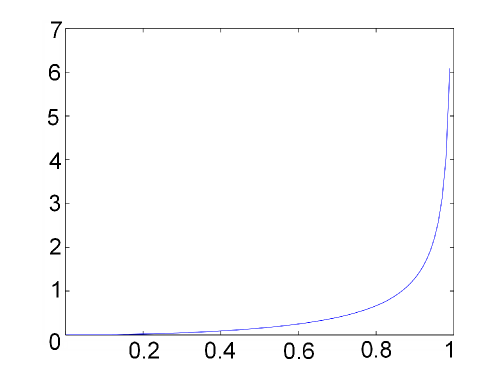}}}
\hfil
\subfigure[]{\raisebox{30.0pt}{\begin{tabular}{|c|c|c|c|c|}
\hline \tiny $\xi$ & \tiny0 & \tiny 1/8 &\tiny 1/7 &\tiny 1/6
\\
\hline \tiny $\Omega(\xi)$ &\tiny 0 &\tiny 0.8\% &\tiny 1\% &\tiny
1.4\%
\\
\hline \tiny $\xi$ & \tiny 1/4 & \tiny 1/3 &\tiny 1/2 &\tiny
$1/\sqrt{2}$
\\
\hline \tiny $\Omega(\xi)$ &\tiny 3\%     &\tiny 6\% &\tiny 15\%
&\tiny  41\%
\\
\hline
\end{tabular}}}
\caption{(a) Graph of $\Omega(\cdot)$; (b) Table of values.}
\label{fig.omega}
\end{figure}

For common bounds $R_j \leq R$, $L_{j,k} = L$, \eqref{basic.disk}
takes the form
\begin{equation*}
\frac{\dis{\bldr(0)}{O_j(0)}}{R} > \Omega(\xi), \quad \frac{L}{R} >
\Upsilon(\xi).
\end{equation*}
Fig.~\ref{fig.ups} illustrates $\Upsilon(\cdot)$.
It follows that if the robot is four times faster than the
obstacles, the spacing between them should exceed nearly negligible
6\% of the obstacle diameter$=$size. If the robot is twice faster,
this percentage increases to $\approx$ the quarter of the size. By
Fig.~\ref{fig.omega}, the requirements to the initial location of
the robot are even more liberal.
\par
\begin{figure}[h]
\centering
\subfigure[]{\scalebox{0.4}{\includegraphics{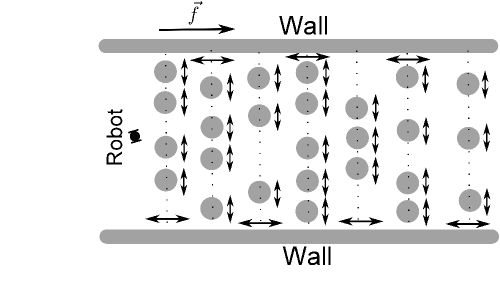}}}
\hfil
\subfigure[]{\scalebox{0.4}{\includegraphics{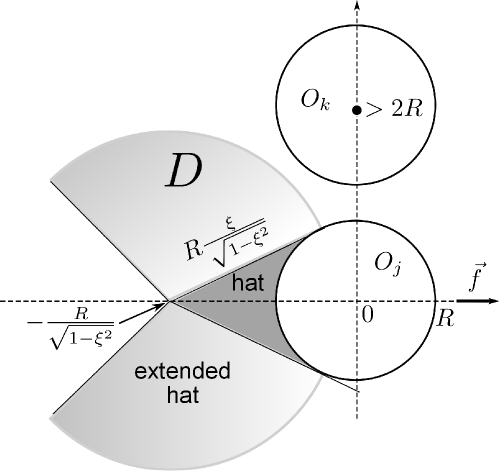}}}
\caption{(a) The corridor; (b) Two close obstacles.} \label{fig.corr}
\end{figure}

{\bf Motion along a corridor obstructed with irregularly and
unpredictably moving disk-shaped obstacles with a common velocity
bound} $v_o$; see Fig.~\ref{fig.corr}(a).
For simplicity only, we assume obstacles with a common radius $R$.
They are partitioned in groups so that the disk centers are always
\textquoteleft vertically\textquoteright (in Fig.~\ref{fig.corr}(a))
aligned within any group. A group may be \textquoteleft
horizontally\textquoteright displaced but its \textquoteleft
horizontal\textquoteright projection does not leave a steady
interval, disjoint with the intervals of other groups. The walls of
corridor are cigar-shaped obstacles (rectangles smoothly
concatenated with two end-disks of proper radii) parallel to
$\vec{f}$. The obstacles move without collisions. Initially the
robot is inside the \textquoteleft rectangular\textquoteright part
of the corridor, to the left of all intervals, and in touch with no
obstacle. The velocity ratio $\xi:= v_o/v<1$, which guarantees
collision avoidance by Theorem~\ref{th.main}.
\par
Starting with the case of one group, we are interested in the
velocity ratio $\xi$ for which
\par
{\bf p)} the robot driven by the properly tuned controller does pass
the corridor in the right direction, no matter how small the minimal
\textquoteleft vertical\textquoteright spacing between the obstacles
is.
\par
Since the hats of the walls lie to the left of their rectangular
parts, a) and b) in Theorem~\ref{th.reach} are fulfilled for walls.
For the disks, a) is true if the extended hat from
Fig.~\ref{fig.fac3}(b) lies in the horizontal strip spanned by the
disk. This holds iff $\frac{\xi}{\sqrt{1-\xi^2}}< 1 \Leftrightarrow
\xi < 1/\sqrt{2} \Leftrightarrow v > \sqrt{2} v_o \approx 1.41 v_o$,
i.e., the robot should be $\approx$ 41\% faster than the obstacles.
Then {\bf p)} holds provided that initially the robot is at a
distance $d_{\text{in}}>R \Omega(1/\sqrt{2}) = \sqrt{2} R \approx
1.41 R$ from the leftmost interval. This extends on many groups if
the distance $d_i$ between the $i$th and $(i+1)$th intervals exceeds
$R \Omega(1/\sqrt{2})$ for any i. We remark that controller tuning
depends only on $d_{\text{in}}$ and $d_i$: it suffices that
$\Delta\left[ \sqrt{2} R \right] = \Delta[0] \in \left( \pi/4,
\arcsin \xi^\star \right)$, where $\xi^\star > 1/\sqrt{2}$ is chosen
so that $d_{\text{in}}, d_i>R \Omega(\xi^\star) \; \forall i$.
\par
In the last example, we omited specifications of the controller
tuning \eqref{specif.cont} for the sake of brevity. They are purely
technical, elementary, and basically in the vein of \eqref{eq.vein}.
\subsection{Thin cigar-shaped obstacles}
\begin{figure}[h]
\centering
\subfigure[]{\scalebox{0.4}{\includegraphics{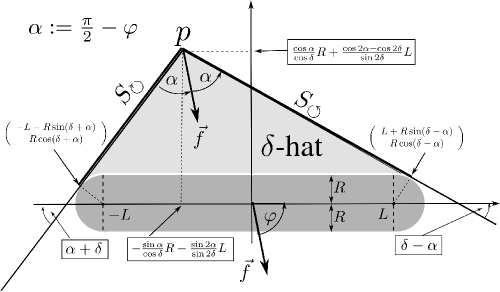}}}
\hfil
\subfigure[]{\scalebox{0.3}{\includegraphics{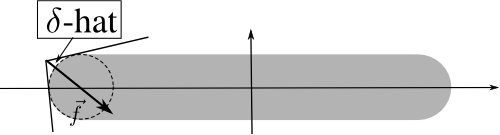}}}
\caption{Cigar and its $\delta$-hat (a) $\alpha:= \frac{\pi}{2} - \varphi \in [-\delta, +\delta]$;
(b) $\alpha \in (\delta, \pi/2]$.}
 \label{fig.cigar}
\end{figure}

{\it Cigar} is the body obtained by smooth concatenation of a
rectangle with two end-disks of proper radii; see
Fig.~\ref{fig.cigar}(a).
Let $\varphi \in (-\pi,\pi)$ stand for the angle from the desired
direction $\vec{f}$ to the cigar centerline. If $|\varphi| \leq
\frac{\pi}{2} - \delta$ or $|\pi - \varphi| \leq \frac{\pi}{2} -
\delta$, the $\delta$-hat and extended hat of the cigar are
identical to those of the respective end-disk; see
Figs.~\ref{fig.cigar}(b) and \ref{fig.fac3}. For $\varphi \in
(\pi/2-\delta, \pi/2+\delta)$, the $\delta$-hat can be found via
elementary geometrical considerations and is shown in
Fig.~\ref{fig.cigar}(a) in lighter dark. The $\delta_1$-extended
$\delta$-hat is obtained by adding the disk sectors swept by the
segments $S_{\circlearrowleft}$ and $S_{\circlearrowright}$ from
Fig.~\ref{fig.cigar}(a) when rotating counter clockwise and
clockwise, respectively, around the hat's vertex $p$ through the
angle $\delta+\delta_1$. If $-\varphi \in (\pi/2-\delta,
\pi/2+\delta)$, the picture is symmetric.
\par
To simplify formulas, we confine ourselves to asymptotic analysis as
$R \to 0$ and consider thin cigars with $R \approx 0$, called {\it
segments}. The subsequent characterization of segment hats is true
with as high precision as desired and the final claims are valid if
$R$ is small enough. Transition to \textquoteleft
non-asymptotic\textquoteright exact claims results from elementary
substitution of the true hats in place of their asymptotic
counterparts.
\par
\begin{figure}[h]
\centering
\scalebox{0.4}{\includegraphics{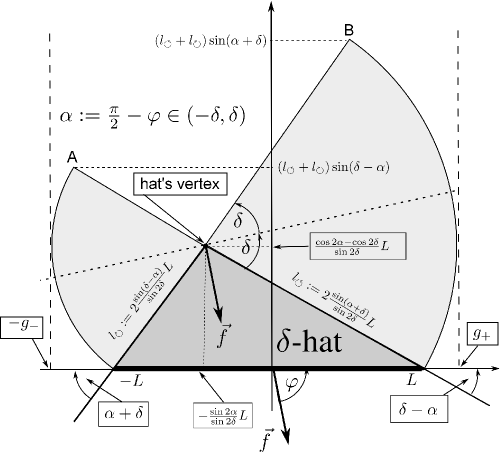}}
\caption{The $\delta$-hat and $\delta$-extended $\delta$-hat of a segment.}
 \label{fig.segment}
\end{figure}

For $|\varphi-\pi/2| < \delta$, the hats of the segment are depicted
in Fig.~\ref{fig.segment} and obtained via elementary geometrical
considerations; the extended hat is composed of the hat and two
slightly shadowed disk sectors.
For $|\varphi+\pi/2| < \delta$, the picture is symmetric. For the
other $\varphi$, the hats are empty.
\par
\begin{figure}[h]
\centering
\subfigure[]{\scalebox{0.4}{\includegraphics{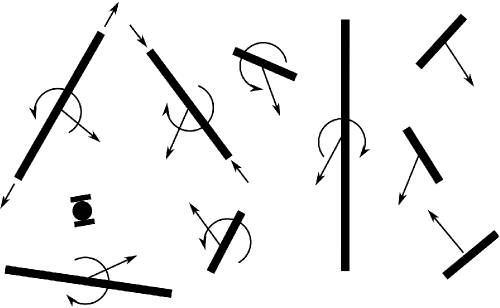}}}
\hfil
\subfigure[]{\scalebox{0.4}{\includegraphics{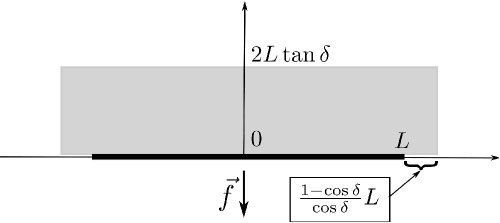}}}
\caption{(a) Scene cluttered with dynamic segments; (b) Rectangle containing the extended hat.}
 \label{fig.segments}
\end{figure}

Now we consider a scene cluttered with moving disjoint segments; see
Fig.~\ref{fig.segments}(a). The $j$th of them has time-varying
half-length $L_j(t)$ and rotates with the angular velocity
$\omega_j(t)$; the velocity of its center is $\vec{v}_j(t)$. For the
sake of brevity, we assume common bounds for the parameters of all
segments:
$$
L_j(t) \leq L, \; |\dot{L}_j(t)| \leq L_1, \; |\omega_j(t)| \leq \omega, \; |\vec{v}_j(t)| \leq v_o \quad \forall t,
$$
and consider controller with only one class: $\Delta_i(\cdot) =
\Delta(\cdot)$.
\par
The speed of any segment's point does not exceed $v_0 + L_1 + L
\omega$. So Theorems~\ref{th.main} and \ref{th.reach} imply the
following.
\begin{prop}
\label{prop.segment} Let the robot be faster than the obstacles:
\begin{equation}
\label{faster.segment}
v_{\Sigma}:= v_0 + L_1 + L
\omega < v.
\end{equation}
If the \textquoteleft parameter\textquoteright $\Delta(\cdot)$ of
the controller is chosen so that $\Delta(0) > \overline{\delta}:=
\arcsin \frac{v_{\Sigma}}{v}$, no collision with the obstacles
occurs. Suppose also that the $\delta$-extended $\delta$-hat of any
segment $j$ (see Fig.~{\rm \ref{fig.segment}}) is always disjoint
with the other segments, and the robot is initially outside their
$\delta$-hats. Moreover let these be true for some $\delta^\star >
\overline{\delta}$. Then the controller can be tuned so that the
robot always drifts in the right direction, i.e., \eqref{drift} is
true. For this to hold, it suffices to pick $\Delta(\cdot)$ so that
\begin{equation}
\label{specif.contseg}
\Delta_i\left[ \frac{1-\cos 2 \delta^\star}{\sin 2 \delta^\star} \right] = \Delta_i[0] \in
(\overline{\delta},\delta^\star).
\end{equation}
\end{prop}
\par
Now we apply this proposition to more special scenarios and, by
sacrificing a part of its content, provide convergence criteria in
terms that are simpler than hats.
\par
{\bf Steady-size segments irregularly and unpredictably moving so
that they remain perpendicular to the desired direction} $\vec{f}$.
The directions of their velocities and so the paths of the centers
are not anyhow restricted. In this case, $L_1=0, \omega=0$ and
\eqref{faster.segment} shapes into $\xi:= v_o/v <1$. It is easy to
see that the $\delta$-extended $\delta$-hat is contained in the
rectangle from Fig.~\ref{fig.segments}(b), where $\tan \delta =
\Gamma(\xi):= \frac{\xi}{\sqrt{1-\xi^2}}$ and
$\frac{1}{2}\frac{1-\cos \delta}{\cos \delta} =
\Xi(\xi):=\frac{1-\sqrt{1-\xi^2}}{2 \sqrt{1-\xi^2}}$ for $\delta =
\arcsin \xi$. Let the spacing between any two obstacles be always no
less than a constant $d_{f}$ in the direction of $\vec{f}$ and no
less than $d_{f}^\bot$ in the perpendicular direction. Then we
arrive at the following.
\begin{claim}
The proposed controller guarantees obstacle avoidance and constant
drift in the right direction \eqref{drift} if there is enough space
between the obstacles $d_f/(2L) > \Gamma(\xi), d^\bot_f/(2L)
> \Xi(\xi)$, the initial distance from the robot to any of them is
no less than $\Gamma(\xi)L$, and $\Delta(\cdot)$ is chosen so that $
\Delta\left[ 1/2 \Gamma(\xi) L + \ve \right] = \Delta[0] \in \left(
\arcsin \xi, \arcsin \xi^\star \right). $ Here $\ve>0$ may be
arbitrarily small and $\xi^\star
> \xi$ is picked so that $d_f/(2L) > \Xi(\xi^\star), d^\bot_f/(2L)
> \Gamma(\xi^\star)$.
\end{claim}
\begin{figure}[h]
\renewcommand{\arraystretch}{1.2}
\centering
\subfigure[]{\scalebox{0.3}{\includegraphics{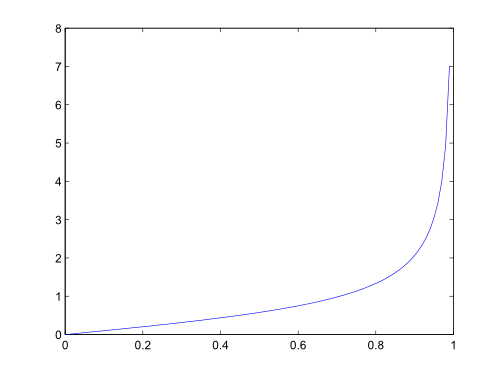}}}
\hfil
\subfigure[]{\raisebox{30.0pt}{\begin{tabular}{|c|c|c|c|c|}
\hline \tiny $\xi$ & \tiny0 & \tiny 1/8 &\tiny 1/7 &\tiny 1/6
\\
\hline \tiny $\Gamma(\xi)$ &\tiny 0 &\tiny 13\% &\tiny 14\% &\tiny
17\%
\\
\hline \tiny $\xi$ & \tiny 1/4 & \tiny 1/3 &\tiny 1/2 &\tiny
$1/\sqrt{2}$
\\
\hline \tiny $\Gamma(\xi)$ &\tiny 26\% &\tiny 35\% &\tiny 58\%
&\tiny  100\%
\\
\hline
\end{tabular}}}
\caption{(a) Graph of $\Gamma(\cdot)$; (b) Table of values.}
\label{fig.gamma}
\end{figure}
\par
\begin{figure}[h]
\renewcommand{\arraystretch}{1.2}
\centering
\subfigure[]{\scalebox{0.3}{\includegraphics{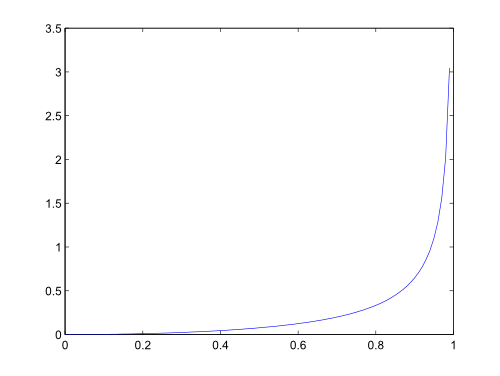}}}
\hfil
\subfigure[]{\raisebox{30.0pt}{\begin{tabular}{|c|c|c|c|c|}
\hline \tiny $\xi$ & \tiny0 & \tiny 1/8 &\tiny 1/7 &\tiny 1/6
\\
\hline \tiny $\Xi(\xi)$ &\tiny 0 &\tiny 0.4\% &\tiny 0.5\% &\tiny
0.7\%
\\
\hline \tiny $\xi$ & \tiny 1/4 & \tiny 1/3 &\tiny 1/2 &\tiny
$1/\sqrt{2}$
\\
\hline \tiny $\Xi(\xi)$ &\tiny 1.5\% &\tiny 3\% &\tiny 8\% &\tiny
20\%
\\
\hline

\end{tabular}}}
\caption{(a) Graph of $\Xi(\cdot)$; (b) Table of values.}
\label{fig.xi}
\end{figure}

The function $\Gamma(\cdot)$ is illustrated in Fig.~\ref{fig.gamma}.
For example, if the robot is at least eight times faster than the
obstacles, spacing between them in the direction of $\vec{f}$ should
exceed only 13\% of the obstacle length for the controller to
succeed. If the robot is twice faster, this percentage increases to
$\approx$ 58\% of the length. As follows from Fig.~\ref{fig.xi}, the
requirements to spacing in the perpendicular direction are much more
liberal.
\par
It is instructive to compare the last claims with the counterexample
from subsect.~\ref{subsect.counter}.
\par
\begin{figure}[h]
\centering
\subfigure[]{\scalebox{0.4}{\includegraphics{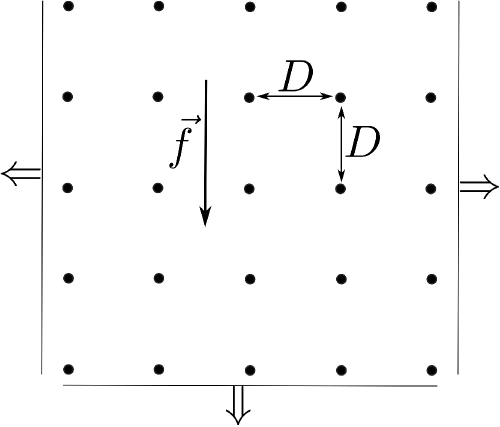}}}
\hfil
\subfigure[]{\scalebox{0.3}{\includegraphics{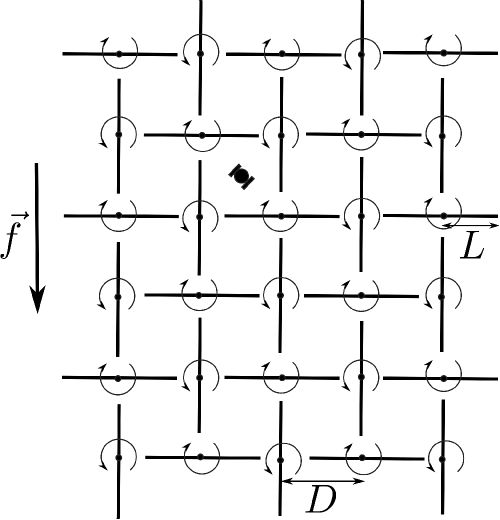}}}
\caption{(a) Cartesian grid; (b) Field of rotating segments.}
\label{fig.grid0}
\end{figure}
{\bf Navigation in the field of rotating segments.} We consider an
infinite Cartesian $D$-grid in the plane; see
Fig.~\ref{fig.grid0}(a). It is not limited along the abscissa axis
and in the negative direction of the ordinate axis, which is
identical to the desired direction $\vec{f}$. At the same time, the
grid has the \textquoteleft upper\textquoteright row.
\par
\begin{figure}[h]
\centering
\subfigure[]{\scalebox{0.4}{\includegraphics{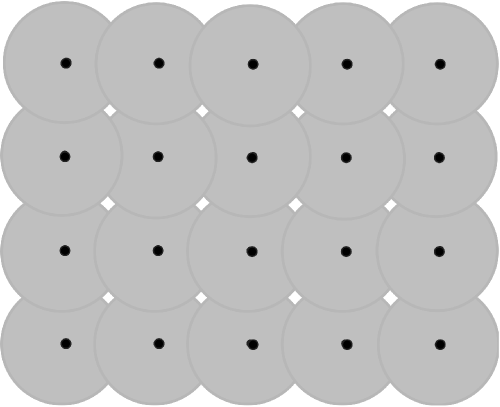}}}
\hfil
\subfigure[]{\scalebox{0.6}{\includegraphics{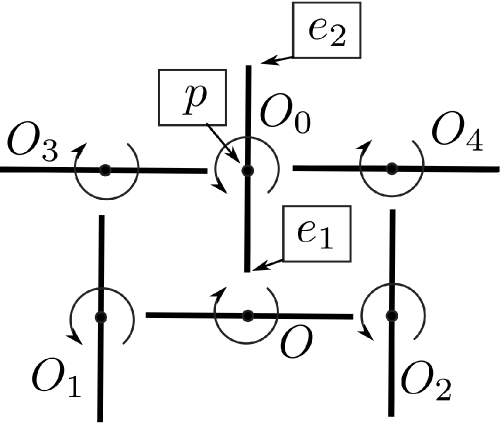}}}
\caption{(a) White areas where collisions with the obstacles are impossible; (b) Five neighboring obstacles.}
\label{fig.grid}
\end{figure}

Every vertex of the grid is a pivot point for a segment with length
$2L$; see Fig.~\ref{fig.grid0}(b). The segments rotate about their
centers with the common angular speed $\omega>0$ in the directions
and from initial orientations shown in Fig.~\ref{fig.grid0}(b). Such
motion of the ensemble is possible if the segments at the ends of
the cell's diagonal do not collide, which is equivalent to
$D>\sqrt{2}L$. To make the things interesting, the focus is on the
case where the paths of the segments overlap $D < 2L$, with the path
being the disk swept by the segment. So the \textquoteleft
path-free\textquoteright part of the plane is composed of infinitely
many compact disconnected areas, shown in Fig.~\ref{fig.grid}(a) in
white. Hence motion through this part with constant drift in the
direction of $\vec{f}$ is not feasible. Initially the robot is in
touch with no segment and above the upper row of the grid.
\par
We assume that the robot is able to overtake any obstacle point
$v>\omega L$. Then by Theorem~\ref{th.main}, the proposed controller
ensures collision avoidance if $\Delta(0)
> \arcsin (L \omega/v)$.
\begin{claim}
\label{claim.rot} The proposed controller provides not only
collision avoidance but also constant drift in the right direction
\eqref{drift} if
\begin{equation}
\label{cond1}
\frac{D}{L} > 2 \sin \delta + \frac{1}{\cos \delta}, \quad
\dis{\bldr(0)}{\mathscr{L}} > \frac{1-\cos 2 \delta}{\sin 2 \delta}L
\end{equation}
for $\delta = \overline{\delta}:= \arcsin(L \omega/v)$, where
$\mathscr{L}$ is the line spanned by the upper row of the grid. For
these to hold, $\Delta(\cdot)$ should be chosen so that $
\Delta\left[ \frac{1-\cos 2 \delta^\star}{\sin 2 \delta^\star}L +
\ve \right] = \Delta[0] \in \left( \overline{\delta}, \delta^\star
\right) $ for some $\ve>0$ and $\delta^\star > \overline{\delta}$
such that \eqref{cond1} is true with $\delta:= \delta^\star$.
\end{claim}
The proof of this claim is given in Appendix~\ref{app.3}.
\par
If $D<2L$, \eqref{cond1} $\Rightarrow \delta<26.26^\circ$ and
$v/(L\omega)>2.26$, i.e., the robot should be $\approx 126\%$ faster
than the obstacles for the constant drift in the right direction.
\par
\begin{figure}[h]
\centering
\subfigure[]{\scalebox{0.3}{\includegraphics{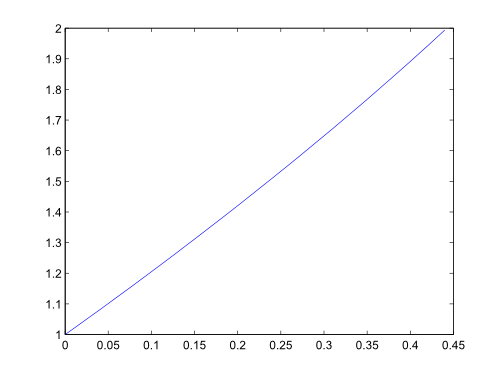}}}
\hfil
\subfigure[]{\raisebox{40.0pt}{\begin{tabular}{|c|c|c|c|}
\hline
\scriptsize $\xi$ & \scriptsize 0 & \scriptsize 1/8 &\scriptsize 1/7
\\
\hline
\scriptsize function &\scriptsize 1 &\scriptsize 26\% &\scriptsize 30\%
\\
\hline
\scriptsize $\xi$ & \scriptsize 1/6 & \scriptsize 1/4 &\scriptsize 1/3
\\
\hline
\scriptsize function &\scriptsize 35\% &\scriptsize 53\% &\scriptsize 73\%
\\
\hline
\end{tabular}}}
\caption{The right-hand side of \eqref{cond1} (a) graph, (b) table of values}
\label{fig.rott}
\end{figure}

The right-hand side of the first inequality from \eqref{cond1}
treated as a function of $\xi=(\omega L)/v$ is illustrated in
Fig.~\ref{fig.rott}.
Since $D/L > \sqrt{2} \approx 1.4142$ (excess $\approx 41\%$), that
inequality is fulfilled for robots that are at least $5.08$ times
faster than the obstacles. For them, the controller ensures constant
drift in the right direction irrespective of how small the white
cells in Fig.~\ref{fig.grid}(a) are. For slower robots, the ratio
$D/L$ should exceed the minimum $\approx 100+41\%$ as is shown in
Fig.~\ref{fig.rott}(b).

\section{Simulations}
\label{sec.sim}
For the velocity controlled vehicle, the control was updated at a period
of $0.1 s$, and 40 rays in Fig.~\ref{fig.scene}(a) were used to detect obstacles. Their edges were associated with the difference $\geq 2 m$ in two
sequentially detected distances. Following a given azimuth was modeled as reaching a faraway target point. The function $\Delta_i(\cdot)$ in \eqref{eq.exfac} was common for all obstacles and specified as a continuous piecewise linear function,
whose fractures are shown in Table~\ref{papr}

\begin{table}
\centering
\begin{tabular}{|c|c|c|c| c|c|c| c| c|}
\hline
\scriptsize $d$(m) & \scriptsize 0 & \scriptsize 0.5 &\scriptsize 1.0 &\scriptsize 1.5 &\scriptsize 2.0 &\scriptsize 2.5 &\scriptsize 3.0 &\scriptsize 100.0
\\
\hline
\scriptsize $\Delta$(rad) & \scriptsize 1.52 & \scriptsize 1.27 &\scriptsize 1.21 &\scriptsize 0.43 &\scriptsize 0.2 &\scriptsize 0.02 &\scriptsize 0.01 &\scriptsize 0.003
\\
\hline
\end{tabular}
\label{papr}
\caption{Parameterization of $\Delta$}
\end{table}

\begin{figure}
\centering
\subfigure[]{\includegraphics[width=0.45\columnwidth]{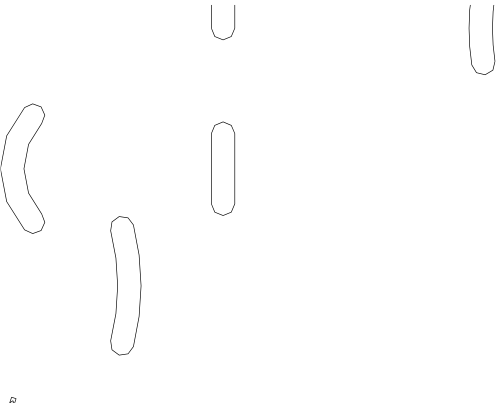}}
\subfigure[]{\includegraphics[width=0.45\columnwidth]{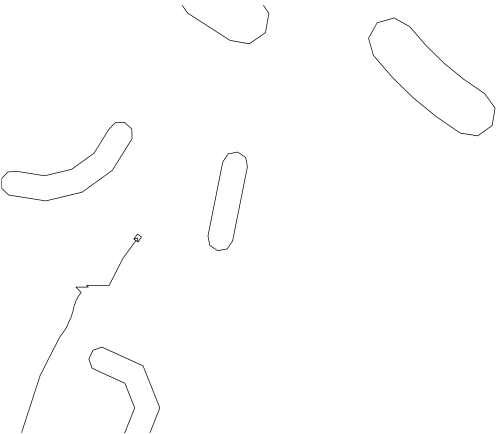}}
\subfigure[]{\includegraphics[width=0.45\columnwidth]{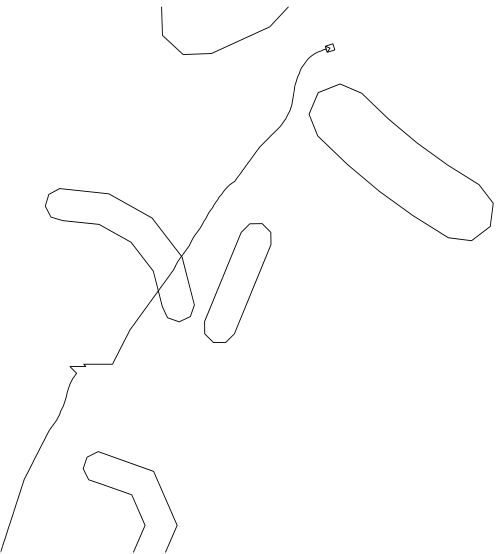}}
\caption{Simulations in a complex environment}
\label{fig:sim2}
\end{figure}

They result from optimization by a genetic algorithm
(MATLAB 2013a, Optimization Toolbox v6.3),
minimizing the time of target reaching while respecting the safety margin $1 m$ to the obstacles in the scenario from Fig.~\ref{fig:sim2}.
\par
\begin{figure}
\centering
\subfigure[]{\includegraphics[width=0.45\columnwidth]{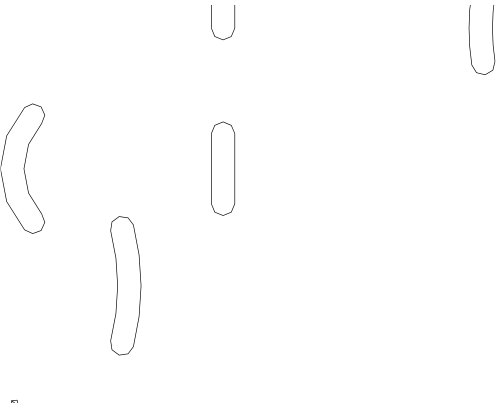}}
\subfigure[]{\includegraphics[width=0.45\columnwidth]{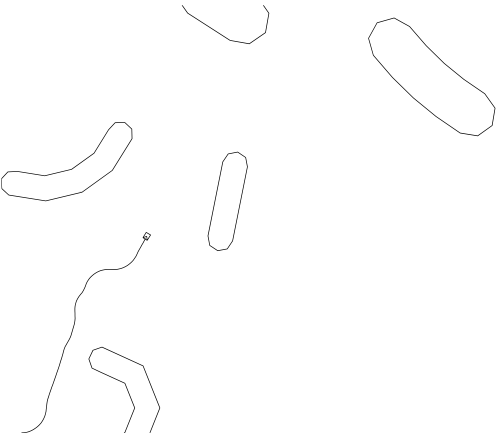}}
\subfigure[]{\includegraphics[width=0.45\columnwidth]{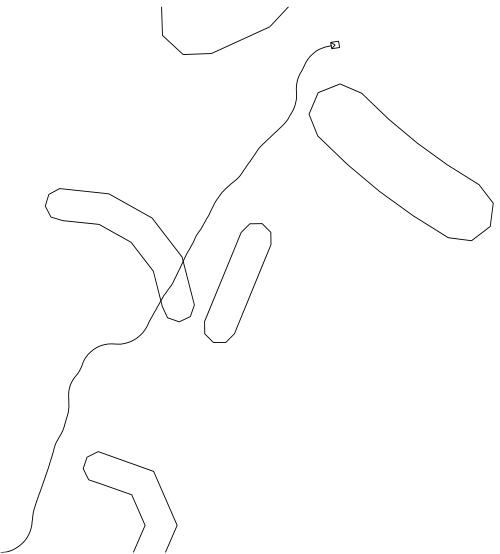}}
\caption{Simulation with vehicle dynamics and noise.}
\label{fig:sim3}
\end{figure}

In the simulation from Fig.~\ref{fig:sim2}, the robot efficiently converges to the target
in a complex scene with many translating and rotating obstacles.
In the simulation from Fig.~\ref{fig:sim3}, the mission was troubled by vehicle dynamics and noise. Unicycle dynamics were used with the turning rate limit of $1 rad s^{-1}$. A random Gaussian disturbance with the standard deviation of $0.1 rad s^{-1}$ was
added to the control input. These entailed no problems for the method, displaying a promising potential to cope with noise and dynamic constraints.

\begin{figure}
\centering
\subfigure[]{\includegraphics[width=0.45\columnwidth]{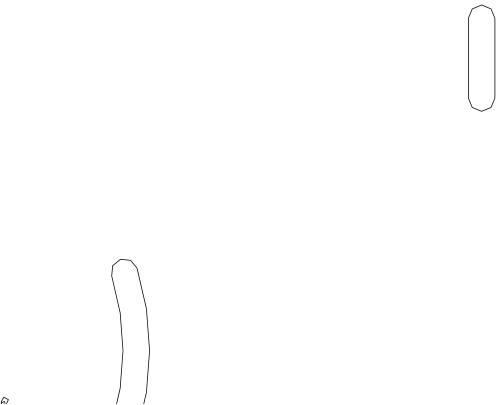}}
\subfigure[]{\includegraphics[width=0.45\columnwidth]{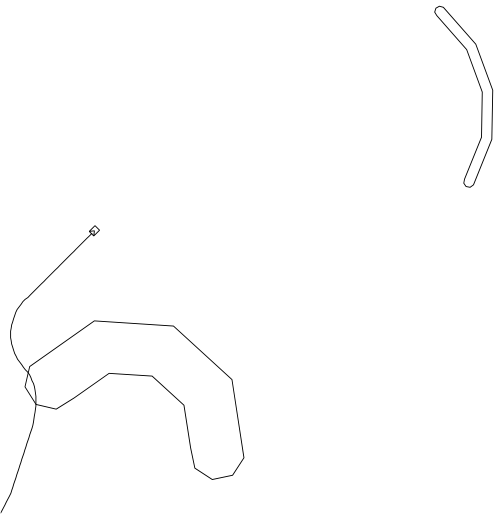}}
\subfigure[]{\includegraphics[width=0.45\columnwidth]{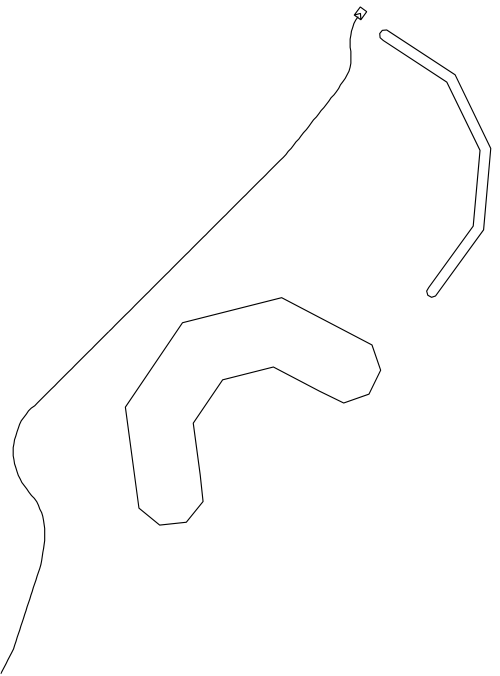}}
\caption{Simulations in a simple environment}
\label{fig:sim1}
\end{figure}

\begin{figure}
\centering
\includegraphics[width=0.9\columnwidth]{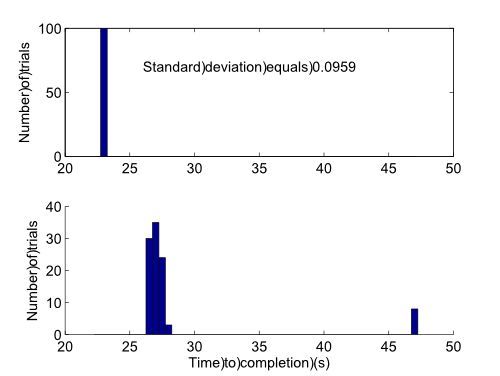}
\caption{Comparison between the PCL and VOM.}
\label{fig:histcomp}
\end{figure}

In the scenario from Fig.~\ref{fig:sim2}, the proposed control law (PCL) was compared with the popular
Velocity Obstacle method (VOM). The latter was given advantage over PCL by access to the obstacle full velocity and knowledge of its future movement within the next $10s$.
VOM also took into account the desired offset from the obstacle and a weighted cost
to compromise progression to the target and separation from obstacles. These parameters were optimized, like for PCL. A Gaussian noise with the standard deviation of $0.1 rad s^{-1}$ was
added to the heading direction. The performance measure was the time
taken to reach the target. By Fig.~\ref{fig:histcomp}, PCL on average outperforms VOM:
VOM took a longer (occasionally, a much longer) path around the obstacles.

\section{Summary}

In this paper, a method for guidance of a Dubins-like
vehicle towards a target in an environment with moving
obstacles is considered. The vehicle is provided with the
relative bearing angle of the target, and the distance to the
nearest point of the obstacle set if it is within the given sensor
range. The proposed control law is composed by biologically
inspired reflex-level rules. Mathematically rigorous analysis
of this law is provided; its convergence and performance are
confirmed by computer simulations.

\bibliographystyle{plain}
\bibliography{Hamidref}

\appendix
\small
\section{Proofs of Lemma~\ref{nec.lem} and Theorem~\ref{th.main}}
\label{app.2} \setcounter{thm}{0}
\renewcommand{\thethm}{A.\arabic{thm}}
{\bf PROOF OF LEMMA~\ref{nec.lem}.} Suppose to the contrary that
$\left[V_{j}^N(\bldr_\ast,t_\ast)\right]_-
> v$ for some $t_\ast$ and $\bldr_\ast \in \partial O_j(t_\ast)$. By
continuity, there exist $\ve
>0$ and $\varkappa >0$ such that
$$
\left[V_{j}^N(\bldr,t)\right]_-
\geq  v +\ve \; \text{if} \; \bldr \in \partial O_j(t), |t-t_\ast| < \varkappa, |\bldr-\bldr_\ast| < \varkappa.
$$
Let the robot start at time $t_\ast$ on the outer normal
$\bldr(t_\ast) = \bldr_\ast - y \vec{N}(\bldr_\ast,t_\ast)$, where
$y >0$ is so small that $\bldr(t_\ast) \not\in O_j(t_\ast)$ and
$|\bldr(t_\ast) - \bldr_\ast| < \varkappa/2$. By
Definition~\ref{def.av}, there exists an admissible velocity profile
$|\vec{v}(t)| \leq v, t \geq t_\ast$ under which
\begin{equation}
\label{vsegda}
d(t) := \dis{\bldr(t)}{O_j(t)} >0 \quad \forall t \geq t_\ast.
\end{equation}
Let $\bldr_\ast(t)$ be the point of $\partial O_j(t)$ closest to
$\bldr(t)$. By Lemma~A.1 \cite{MaChaSa12a} for $t \in [t_\ast,
t_\ast+\tau]$ and $y \approx 0, \tau \approx 0$,
$$
|\bldr_\ast(t)-\bldr_\ast| < \varkappa, \quad |t-t_\ast| < \varkappa \Rightarrow
\left[V_{j}^N[\bldr_\ast(t),t]\right]_- \geq  v +\ve.
$$
At the same time according to (A.1) in \cite{MaChaSa12a},
$$
\dot{d}(t) = V^N_j(t) - v^N(t)  \leq - (v+\ve) +v = - \ve.
$$
This and \eqref{vsegda} yield $ \ve \tau \leq d(0) = y $, though $y$
can be chosen as close to $0$ as desired. This contradiction proves
the lemma. \epf
\par
\begin{figure}
\subfigure[]{
\scalebox{0.25}{\includegraphics{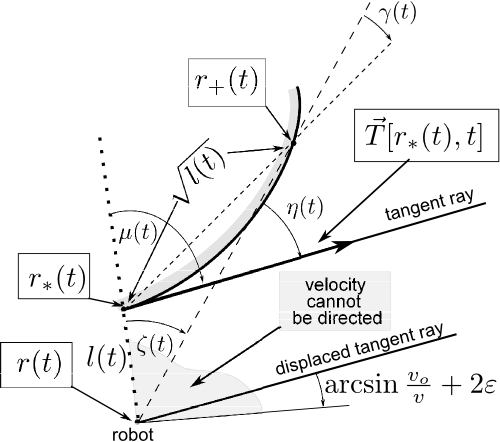}}}
\hfill
\subfigure[]{
\scalebox{0.25}{\includegraphics{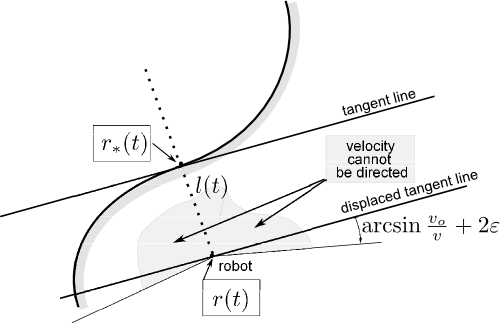}}}
\caption{Angles}
\label{fig.angles}
\end{figure}

{\bf PROOF OF THEOREM~\ref{th.main}.} Suppose the contrary: $G:=\{t
\geq 0: \bldr(t) \not\in \bigcup_i O_i(t)\} \neq [0,\infty)$. Since
$G$ is open in $[0,\infty)$, its leftmost connected component has
the form $[0,\tau)$, where $\bldr(\tau) \in O_k(\tau)$ for some $k$.
Let $\bldr_\ast(t)$ be the point of $\partial O_k(t)$ closest to
$\bldr(t)$ for $t\leq \tau$, and $l(t):= |\bldr(t)-\bldr_\ast(t)|$.
For $t \in (0,\tau)$, we introduce (see Fig.~\ref{fig.angles}(a))
the point $\bldr_\pm(t) \in \partial O_k(t)$ that is
counter-clockwise/clockwise ahead of $\bldr_\ast(t)$ and separated
by $\sqrt{l(t)}$ along the boundary $\partial O_k(t)$, the polar
angle  $\alpha(t)$ of $\overrightarrow{\bldr(t), \bldr_\ast(t)}$ in
the robot's frame and the angle $\zeta(t)$ from
$\overrightarrow{\bldr(t), \bldr_\ast(t)}$ to
$\overrightarrow{\bldr(t), \bldr_+(t)}$, the angle $\gamma(t)$ from
$\overrightarrow{\bldr(t), \bldr_+(t)}$ to
$\overrightarrow{\bldr_\ast(t), \bldr_+(t)}$, the angle $\eta(t)$
from $\overrightarrow{\bldr_\ast(t), \bldr_+(t)}$ to
$\vec{T}[\bldr_\ast(t),t]$, and the angle $\mu(t)$ from
$\overrightarrow{\bldr(t), \bldr_\ast(t)}$ to
$\vec{T}[\bldr_\ast(t),t]$.
Here $\mu(t) = \zeta(t)+\gamma(t)+\eta(t)$ and $\eta(t) \to 0,
\gamma(t) \to 0$ as $t \to \tau-$ since the vector $\vec{T}$ is
tangent and $|\bldr_\ast(t) - \bldr_+(t)|/|\bldr(t) - \bldr_\ast(t)|
\sim \sqrt{l(t)}/ l(t) \to \infty $ as $t \to \tau-$. Thus $ \mu(t)
- \zeta(t) \to 0 \quad \text{as} \quad t \to \tau- $. By
\eqref{safety.req}, there exist $\ve>0, \varkappa>0$ such that
$\Delta_i(d) - 2 \ve \geq
\overline{\delta}:=\arcsin\frac{v_o^i}{v}\; \forall d \leq
\varkappa$, where $C_i\ni k$. Now we focus on $t<\tau$ so close to
$\tau$ that $|\mu(t) - \zeta(t)| \leq \ve$ and $l(t) \leq
\varkappa$, which implies $d_k^{\,\min}(t) \leq \varkappa$ in
\eqref{eq.exfac}.
\par
The robot's view in the directions $\alpha \in [\alpha(t)-\zeta(t),
\alpha(t)]$ is obstructed by $O_k(t)$. Hence
$I(t):=[\alpha(t)-\zeta(t)-\Delta_i[d_k^{\,\min}(t)], \alpha(t)]$ is
covered by the range $\widehat{A}_k(t)$ of the extended facet
$\widehat{F}_k$. The respective part of $\widehat{F}_k$ cannot be
shadowed by other obstacles for $t \approx \tau$. Indeed as $t \to
\tau-$, the locus of \textquoteleft shadowing\textquoteright points
collapses into a part of $\partial O_k(\tau)$ and so the distance
from the locus to the other obstacles is lower limited by a non-zero
value since the obstacles do not collide. Hence for the polar angle
$\alpha_v(t)$ of the robot's velocity $\vec{v}(t)$, we have
$\alpha_v(t) \not\in I(t) \supset
[\alpha(t)-\zeta(t)-\overline{\delta}-2 \ve, \alpha(t)] \supset
[\alpha(t)-\mu(t)-\overline{\delta}-  \ve, \alpha(t)]$. By replacing
$\bldr_+(t) \mapsto \bldr_-(t)$ in the foregoing, we see that
$\alpha_v(t) \not\in [\alpha(t), \alpha(t)+\mu(t)+\overline{\delta}+
\ve]$.
\par
Overall, $\vec{v}(t)$ is not directed to the shadowed sector from
Fig.~\ref{fig.angles}(b), i.e., the angle subtended by the velocity
$\vec{v}(t)$ and the unit inner normal $N[\bldr_\ast(t), t]$ is no
less than $\pi/2+\arcsin\frac{v_o}{v}+  \ve$. So the related normal
component $v^N(t)$ of $\vec{v}(t)$
\begin{multline*}
v^N(t) \leq - v \sin \left( \arcsin\frac{v_o}{v}+ \ve \right) \leq -
v_o - \ve_\ast, \; \text{where}
\\
\ve_\ast:= v \left[ \sin \left( \arcsin\frac{v_o}{v}+ \ve \right) -
\sin \arcsin\frac{v_o}{v} \right] >0.
\end{multline*}
At the same time according to (A.1) in \cite{MaChaSa12a},
$$
\dot{l}(t) = V^N_j(t) - v^N(t)  \overset{\text{\eqref{spee.req}}}{\geq}
- v_o + v_0 + \ve_\ast = \ve_\ast >0
$$
in violation of $l(t)>0 \; \forall t \in (0,\tau)$ and $l(t) \to 0 $
as $t \to \tau-$. The contradiction obtained completes the proof.
\section{Proof of Theorem~\ref{th.reach}}
\label{app.1}
\setcounter{thm}{0}
\renewcommand{\thethm}{B.\arabic{thm}}
\renewcommand{\theobservation}{B.\arabic{observation}}
We start with a deeper insight into the $\delta$-hats of convex
bodies, where $\delta \in (0,\pi/2)$. We fix $t$ and $j$, focus on
$O:=O_j(t)$, assume the polar angle of $\vec{f}$ to be zero, and add
the adjective {\it free} to signal that the facet is calculated in
the absence of the other obstacles.
\par
\begin{figure}[h]
\centering
\subfigure[]{\scalebox{0.15}{\includegraphics{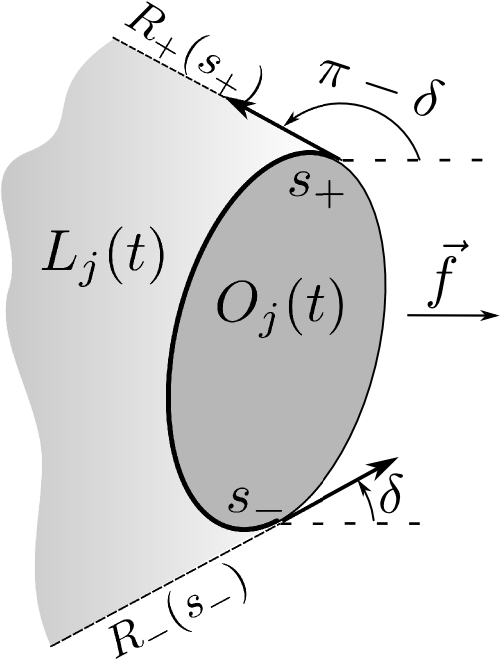}}}
\hfil
\subfigure[]{\scalebox{0.15}{\includegraphics{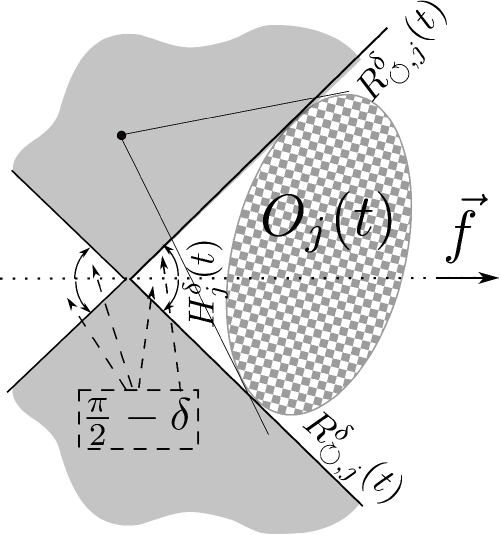}}}
\hfil
\subfigure[]{\scalebox{0.15}{\includegraphics{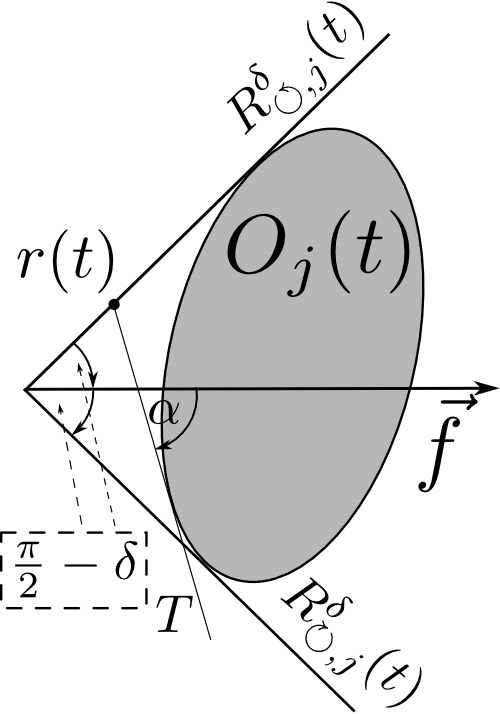}}}

\caption{(a) Points for which the direction of $\vec{f}$
is shadowed by the extended facet; (b) Continuity points; (c) Free extended facet.}
\label{fig.shadow}
\end{figure}

Let $\bldr = \rho(s)$ be a natural parametric representation of the
boundary $\partial O$, where $s$ is the arc length and $O$ is to the
left as $s$ ascends, $\vec{T}(s):= \frac{d \rho}{ds}(s)$ and
$\vec{N}(s)$ be the unit normal to $\partial O$ directed inwards
$O$. The ray $R_\pm(s):= \{\bldr = \bldr_\pm(x,s):= \rho(s) \pm x
\vec{T}(s): x \geq 0\}$ is the locus of points from which the
uppue/lower edge of the free visible facet of $O$ is given by
$\rho(s)$. By the Frenet-Serrat equations $\bldr_\pm^\prime(x,s) =
\vec{T}(s) \pm x \varkappa(s) \vec{N}(s)$, where the signed
curvature $\varkappa(s) \geq 0$ since $O$ is convex. Hence as $s$
ascends, the ray $R_\pm(s)$ (nonstrictly) displaces to the left with
respect to its current position observed from the ray origin. By the
same argument, $\vec{T}(s)$ rotates counter-clockwise, and makes one
full turn as $s$ runs the entire perimeter of $\partial O$. As a
result, we arrive at the following.
\begin{observation}
\label{obs1} {\bf a)} Let $s_-$ and $s_+$ be values of $s$ for which
the polar angle $\alpha(s)$ of $\vec{T}(s)$ equals $\delta$ and $\pi
- \delta$, respectively. The locus $L_j(t)$ of points $\bldr \not\in
O_j(t)$ whose view in the direction of $\vec{f}$ is obstructed by
the free $\delta$-facet of $O_j(t)$ is bounded by $R_+(s_+)$,
$R_-(s_-)$, and the respective part of $\partial O_j(t)$
(see Fig.~{\rm
\ref{fig.shadow}(a)}).
\\
{\bf b)} Let $s_\pm^\delta$ be the largest/smallest value of $s$ for
which $\alpha(s)=3/2\pi\pm \delta$, respectively. The $\delta$-hat
$H^\delta_j(t)$ of $O_j(t)$ is bounded by $R_+(s_-^\delta)$,
$R_-(s_+^\delta)$, and the least part of $\partial O_j(t)$ that
connects these rays and excludes $O_j(t)$ from the hat.
\\
{\bf c)} $H_j^{\delta_1}(t) \subset H_j^{\delta_2}(t),
s_-^{\delta_1} \geq s_-^{\delta_2}, s_+^{\delta_2} \geq
s_+^{\delta_1}$ if $\delta_1 \leq \delta_2$.
\\
{\bf d)} $s_-^{\delta} \to s_-^{\delta_\star}, s_+^{\delta} \to
s_+^{\delta_\star}$ as $\delta \to \delta_\star+$, and the ends of
the segments $S_{\circlearrowright,j}^\delta(t),
S_{\circlearrowleft,j}^\delta(t)$ from Fig.~\ref{fig.hat} converge
to the respective ends of
$S_{\circlearrowright,j}^{\delta_\star}(t),
S_{\circlearrowleft,j}^{\delta_\star}(t)$ uniformly over $t$.
\end{observation}
\begin{corollary}
\label{cor.look} The interiors of the sectors shadowed in Fig.~{\rm
\ref{fig.shadow}(b)} do not contain points for which the obstacle
edges have equal angular discrepancy with respect to $\vec{f}$.
\end{corollary}
This holds by c) in Observation~\ref{obs1} since such points are
vertices of hats. Exactly at these points the control law is
discontinuous at time $t$ if the effect of the other obstacles can
be neglected.
\begin{lem}
\label{lem.maxdis} The maximal distance from a point of the hat
$H_j^\delta(t)$ to $O_j(t)$ equals the height $h_j^\delta(t)$ of the
hat.
\end{lem}
\pf Since $O_j(t)$ is convex, so is the distance function $\bldr
\mapsto \dis{\bldr}{O_j(t)}$. Hence its maximum over the triangle
$\mathscr{T}(t)$ spanned by the vertex $p_j^\delta(t)$ of the angle
$\sphericalangle_j^\delta(t)$ and the points from 1) and 2) on
page~\pageref{semg} is attained a vertex $p \not\in O_j(t)$, which
is $p_j^\delta(t)$. It remains to note that this maximum equals that
over the hat. \epf
\par
The ray $R_{\aleph,j}^\delta(t), \aleph =\circlearrowleft,
\circlearrowright$ evidently moves without rotation. Its {\it normal
velocity} $V_\aleph(t)$ is that with respect to itself; with the
positive direction being to the left when observing from the ray
origin.
\begin{lem}
\label{lem.normmal} Putting $\sigma_{\circlearrowleft}:=-1,
\sigma_{\circlearrowright}:=1$, we have for $\aleph
=\circlearrowleft, \circlearrowright$
$$
V_{\aleph}(t+/-)= \sigma_\aleph \underset{\bldr \in
R_{\aleph,j}^\delta(t) \cap O_j(t)}{\min/\max} V_j^N(\bldr,t).
$$
\end{lem}
\pf Let $\aleph =\circlearrowleft$ and $\rho(s)$ be the natural
parametric representation of the reference configuration boundary
$\partial O_{\ast,j}$. We put $\Phi_j[s,t] := \Phi_j[\rho(s),t],
\vec{V}_j(s,t):= \vec{V}_j[\rho(s),t]$, etc. There are $s_1,s_2$
such that $\Phi_j[s,t] \in R_{\circlearrowleft,j}^\delta(t) \cap
O_j(t) \Leftrightarrow s \in [s_1,s_2]$. Let
$\vec{n}_\circlearrowleft$ be the positively oriented unit vector
normal to $R_{\circlearrowleft,j}^\delta(t)$. The line
$$
\spr{\bldr}{\vec{n}_\circlearrowleft} = \zeta(t):=\max_{\bldr_\ast \in
O_j(t)} \spr{\bldr}{\vec{n}_\circlearrowleft} = \max_{s} \spr{\vec{n}_\circlearrowleft}{\Phi_j(s,t)}
$$
contains $R_{\circlearrowleft,j}^\delta(t)$ and so
$V_{\circlearrowleft}(t\pm) = \dot{\zeta}(t\pm)$. By Corollary~2 in
\cite[Sect.~2.8]{Clarke83} combined with Definition~2.3.4 and
Proposition~2.1.2 \cite{Clarke83}, the one-sided derivatives on the
right do exist and
$$
\dot{\zeta}(t+/-) = \max_\mu/\min_{\mu} \int_{s_1}^{s_2} \frac{\partial \spr{\vec{n}_\circlearrowleft}{\Phi_j(s,t)}}{\partial t} \; \mu(ds),
$$
where $\mu$ ranges over all probability measures on $[s_1,s_2]$. So
this $\max/\min$ are equal to $\max/min_{s \in [s_1,s_2]}$ of the
integrand
$$
\frac{\partial \spr{\vec{n}_\circlearrowleft}{\Phi_j(s,t)}}{\partial t} =
\spr{\vec{n}_\circlearrowleft}{\vec{V}_j(s,t)} \overset{\text{(a)}}{=} - V_j^N(s,t),
$$
where (a) holds since $\vec{n}_\circlearrowleft = - \vec{N}_j(s,t)\;
\forall s \in [s_1,s_2]$. This completes the proof. The case
$\aleph=\circlearrowright$ is considered likewise. \epf
\par
Now we start to examine the robot driven by the proposed control
law. The $\delta$-{\it facet} of the obstacle is given by
\eqref{eq.exfac} and \eqref{eq.exfac1}, where $\delta_k$ is replaced
by $\delta$ everywhere.
\begin{figure}[h]
\centering
\subfigure[]{\scalebox{0.15}{\includegraphics{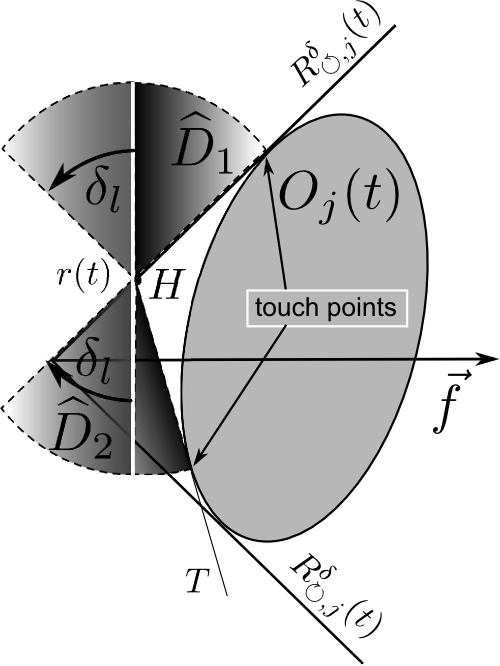}}}
\hfil
 \subfigure[]{\scalebox{0.15}{\includegraphics{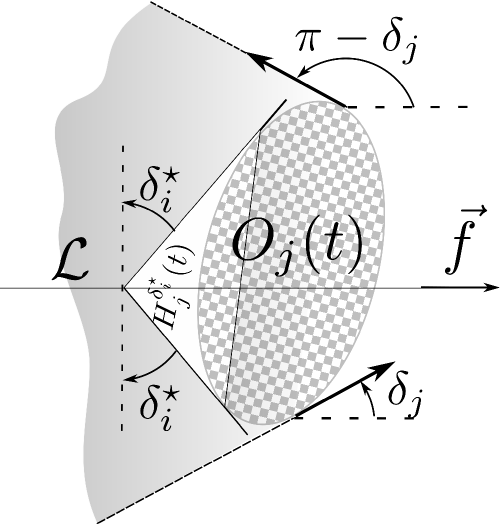}}}
\caption{(a) Disk sectors that should be free; (b) Larger disk sectors; (c) Location of the robot.}
\label{fig.fac}
\end{figure}
\begin{lem}
\label{lem.outside} {\bf a)} If the robot is on
$S_{\circlearrowleft,j}^{\delta_i^\star}(t) \setminus
\{p_j^{\delta_i^\star}(t)\}, j \in C_i$, the normal component of its
relative velocity with respect to each of the rays
$R_{\circlearrowleft,j}^{\delta_i^\star}(t),
R_{\circlearrowright,j}^{\delta^\star_i}(t)$ is no less than $v \sin
\delta^i - v^i_o$, where $\delta^i:= \Delta_i[0]$. {\bf b)} If the
robot is on $S_{\circlearrowright,j}^{\delta_i^\star}(t) \setminus
\{p_j^{\delta_i^\star}(t)\}, j \in C_i$, the above normal component
does not exceed $-(v \sin \delta^i - v^i_o)$.
\end{lem}
\pf We focus on (a), (b) is proved likewise. Since $\bldr(t) \in
S_{\circlearrowleft,j}^{\delta_i^\star}(t)$ and $O_j(t)$ is convex,
$O_j(t)$ obstructs the view in the direction of $\vec{f}$, lies in
the angle $\sphericalangle$ subtended by
$R_{\circlearrowleft,j}^{\delta_i^\star}(t)$ and the tangent ray $T$
from Fig.~\ref{fig.shadow}(c), and its free facet has the form
$\mathscr{A}_j = \left( -\alpha,
\frac{\pi}{2}-\delta_i^\star\right)$. So $\alpha
> \frac{\pi}{2}-\delta_i^\star$, whereas $\alpha \leq
\frac{\pi}{2}+\delta_i^\star$ since $\sphericalangle$ does not
exceed $\pi$ radian. By \eqref{eq.delta}, \eqref{eq.exfac}
\eqref{delta.range} and Lemma~\ref{lem.maxdis}, the range of the
free extended facet $\widehat{\mathscr{A}}_j = \left( -\alpha -
\delta^i, \frac{\pi}{2}-\delta_i^\star+\delta^i\right)$ does not
cover the circle since its angular span $\frac{\pi}{2}
-\delta_i^\star+\delta^i +[\alpha+\delta^i] \leq \pi + 2 \delta^i <
2 \pi$ thanks to $\delta^i < \frac{\pi}{2}$.
So by \eqref{delta.range}, the robot's absolute velocity $\vec{v}$
has the polar angle $\frac{\pi}{2}-\delta_i^\star+\delta^i$ unless a
sector $D_\nu$ from Fig.~\ref{fig.fac}(a) is partly shadowed by an
extended facet of another obstacle $O_k(t)$. Now we show that this
does not hold.
\par
Suppose the contrary. Then obstacle $O_k(t)$ (of class $C_l$)
contains a point $p$ whose rotation through an angle $\leq \delta_l$
about $\bldr(t)$ goes through $\br^2 \setminus O_j(t)$ and ends on
$D_1 \cup D_2$. The locus of such points $p \in \br^2$ lies in the
union of the domain $H$ from Fig.~\ref{fig.fac}(b) and two larger
sectors $\widehat{D}_1$ and $\widehat{D}_2$ that are obtained from
$D_1$ and $D_2$, respectively, by extra rotation of the radius
through angle $\delta_l$. However by \eqref{delta.range}, this union
is a subset of the $\delta_l^\star$-extended $\delta_i^\star$-hat of
$O_j(t)$. Thus $O_k(t)$ intersects this hat, in violation of a) in
Theorem~\ref{th.reach}. Hence $\vec{v}$ has the polar angle
$\frac{\pi}{2}-\delta_i^\star+\delta^i$ indeed.
\par
So $\vec{v}$ subtends the angle $\delta^i$ with the ray
$R_{\circlearrowleft,j}(t)$ and the angle $2 \delta_i^\star -
\delta^i$ with $-R_{\circlearrowright,j}(t)$. So the normal
components of $\vec{v}$ are
$\spr{\vec{v}}{\vec{n}_{\circlearrowleft}} = v \sin \delta^i$ and
$\spr{\vec{v}}{\vec{n}_{\circlearrowright}} = v \sin
(2\delta^\star_i - \delta^i)$, where
$
\sin (2\delta^\star_i - \delta^i) - \sin \delta^i =
2 \sin (\delta^\star_i - \delta^i) \cos \delta^\star_j >0,
$
and so no less than $v \sin \delta^i$ in both cases. It remains to
note that both rays are translated at the speed $\leq v^i_o$ by
\eqref{spee.req1} and Lemma~\ref{lem.normmal}. \epf
\par
By \eqref{delta.ineq}, b) in Theorem~\ref{th.reach}, and
Corollary~\ref{cor.look}, we get the following.
\begin{corollary}
\label{cor.outside} The robot is always outside the
$\delta_i^\star$-hats of all obstacles of class $C_i$.
\end{corollary}
\par
{\bf PROOF OF THEOREM~\ref{th.reach}:} If at time $t$ the ray $R$
emitted from the robot in the direction of $\vec{f}$ hits no
extended facet, the robot's velocity $\vec{v} = v \vec{f}$ and so
\eqref{drift} does hold with $\ve:=v$. Suppose that $R$ hits an
extended facet; let $j \in C_j$ be the index of the first of them
and $\delta_j$ be defined from \eqref{eq.exfac}. By
Corollary~\ref{cor.outside}, the robot is not in the
$\delta_i^\star$-hat of the $j$th obstacle. So in view of a) in
Observation~\ref{obs1}, the robot is located in the shadowed domain
from Fig.~\ref{fig.fac}(c).
\par
Let the robot is above the line $\mathcal{L}$ from
Fig.~\ref{fig.fac}(c). Let $\alpha_j^+$ be the polar angle of the
upper edge of the free facet of $O_j(t)$. Then $- \delta_j \leq
\alpha^+_j \leq \pi/2 - \delta_i^\star\Rightarrow 0 \leq
\alpha^+_j+\delta_j \leq \pi/2 - \delta_i^\star + \delta_j$, where
$\alpha^+_j+\delta_j$ is the upper end of the extended facet. Since
$\delta_i^\star > \delta^i \geq \delta_j$ due to
\eqref{delta.range}, its deviation from the angle $0$ of $\vec{f}$
does not exceed $\beta \leq \pi/2 - \delta_i^\star + \delta^i <
\pi/2$. So were the robot's velocity $\vec{v}$ directed to that
upper end, the angle $\gamma$ between $\vec{v}$ and $\vec{f}$ would
be acute $\gamma \leq \beta$. That direction can be dismissed by
either the lower edge of the same obstacle or by another obstacle.
However this may result only in decrease of $\gamma$. Thus in any
case, $\gamma \leq \beta$ and so
$$
\spr{\vec{v}}{\vec{f}} = v \cos \gamma \geq v \cos \left( \frac{\pi}{2} - \delta_i^\star +
\delta^i \right) = v \sin \left( \delta_i^\star-\delta^i \right) >0,
$$
i.e., \eqref{drift} does hold with $\ve:=v \min_i \sin \left(
\delta_i^\star-\delta^i \right)$.
The case where the robot is below $\mathcal{L}$ is considered likewise. \epf
\section{Proof of Claim~\ref{claim.rot}}
\label{app.3}
\setcounter{thm}{0}
\renewcommand{\thethm}{c.\arabic{thm}}
\renewcommand{\theobservation}{C.\arabic{observation}}
Since the right-hand side of
the second inequality from \eqref{cond1} estimates the height of the
$\delta$-hat from above by Fig.~\ref{fig.segment}, the robot is
initially outside the $\delta^\star$-hats of all horizontal
segments. The same is true for the vertical segments since their
hats are empty.
\par
It remains to check that the extended hat of any obstacle $O$ is
always disjoint with the other obstacles. In doing so, it suffices
to examine only five neighbors $O_0, \ldots,O_4$ of $O$ (see
Fig.~\ref{fig.grid}(b)) for $|\alpha| < \delta$ since otherwise the
hat is empty.
\par
$\mathbf{O_0}.$ The $\delta$-extended $\delta$-hat is disjoint with
$O_0$ if in Fig.~\ref{fig.segment}, the ends $e_i$ of $O_0$ lie
above the lines spanned by the segments of the lengths
$l_\circlearrowleft$ and $l_\circlearrowright$. Analytically this
means that
$$
\begin{array}{r}
-x_i \sin(\delta+\alpha) + y_i \cos(\delta+\alpha) \geq L\sin(\delta+\alpha)
\\
x_i \sin(\delta-\alpha) + y_i \cos(\delta-\alpha) \geq L\sin(\delta-\alpha)
\end{array}
\; \begin{array}{l}
\forall  \alpha \in (-\delta, \delta)
\\
\forall i=1,2
\end{array}.
$$
Here $x_i = -D\sin \alpha - (-1)^iL \sin 2 \alpha, y_i = D \cos
\alpha +(-1)^i L \cos 2 \alpha$ are the coordinates of $e_i$ in the
local frame of $O$ (see Fig.~\ref{fig.segment}). Via elementary
triginometrical identities, these inequalities are shaped into
$$
D \cos \delta \geq L \max_{\alpha\in [0,\delta]}\max\left\{
\begin{array}{l}
\sin(\delta+\alpha) +| \cos(\delta-\alpha)|
\\
\sin(\delta-\alpha) +|\cos(\delta+\alpha)|
\end{array}\right\}  .
$$
This is true by \eqref{cond1} since $\max_{\alpha} \ldots$ does not
exceed $\sin 2 \delta +1$.
\par
$\mathbf{O_1, O_2}.$ We focus on $O_1$, $O_2$ is considered
likewise. The pivot of $O_1$ has the coordinates $x=-D\cos \alpha ,
y=-D\sin\alpha $; the unit vector normal to $O_1$ is $(\cos 2
\alpha, \sin 2 \alpha)$. So the line $\mathscr{L}_1$ spanned by
$O_1$ is given by $x \cos 2 \alpha + y \sin2\alpha = - D \cos
\alpha$. The $\delta$-extended $\delta$-hat from
Fig.~\ref{fig.segment} is disjoint with $O_1$ if the distance from
the vertex $V$ of the $\delta$-hat from Fig.~\ref{fig.segment} to
$\mathscr{L}_1$ exceeds the radius $l_{\circlearrowright}$ of the
left shadowed disk segment. With regard to the coordinates of $V$
given in Fig.~\ref{fig.segment}, this means that for all $\alpha \in
[-\delta,\delta]$,
$$
\begin{array}{c}
-L \frac{\sin 2 \alpha}{\sin 2\delta} \cos 2 \alpha + L \frac{\cos 2 \alpha-\cos 2 \delta}{\sin 2 \delta}\sin 2 \alpha
>
\\
-D \cos \alpha + 2 L \frac{\sin(\delta-\alpha)}{\sin 2 \delta}
\\
\Updownarrow
\\
\frac{D}{L} > \Theta(\alpha):= \frac{2}{\sin 2 \delta} \left[ \cos 2 \delta \sin \alpha +
\sin \delta - \cos \delta \tan \alpha\right] .
\end{array}
$$
By elementary calculus, $\Theta^\prime(\alpha) \leq 0$ and so
$\max_{\alpha \in [-\delta,\delta]} \Theta(\alpha) =
\Theta(-\delta)$, which shapes the condition into
$$
\frac{D}{L} > \frac{2}{\sin 2 \delta} \left[ -\cos 2 \delta \sin
\delta + 2\sin \delta \right] = \frac{3 -2 \cos^2 \delta}{\cos
\delta}
$$
The difference between the expression on the right and the
right-hand side of the first inequality from \eqref{cond1} equals
\begin{multline*}
2 \sin \delta + \frac{1}{\cos \delta} - \frac{3 -2 \cos^2
\delta}{\cos \delta} = \frac{\sin 2\delta + \cos 2\delta- 1 }{\cos
\delta} =
\\
\sqrt{2}\frac{\sin (2\delta + 45^\circ)- \sin 45^\circ }{\cos
\delta} >0 \quad \text{since}\; 0<\delta<26.26^\circ.
\end{multline*}
Thus $O_1$ is disjoint with the $\delta$-extended $\delta$-hat if
\eqref{cond1} holds.
\par
$\mathbf{O_3, O_4}.$ We focus on $O_3$. Since $O_3$ is always
parallel to $O$, it suffices to show that the relative ordinate
$y_p=D(\cos \alpha-\sin\alpha )$ of its pivot exceeds those of A and
B in Fig.~\ref{fig.segment}:
$$
D(\cos \alpha-\sin\alpha ) > L \frac{\cos\alpha}{\cos \delta} \max\{\sin(\delta-\alpha); \sin(\delta+\alpha)\}
$$
for all $\alpha \in [-\delta,\delta]$. Only $\alpha \in [0,\delta]$
can be examined since for them, $\alpha := -\alpha$ keeps the
inequality true. Then it shapes into
$$
D(\cos \alpha-\sin\alpha ) > L \frac{\cos\alpha}{\cos \delta} \sin(\delta+\alpha)
$$
Since the left-and right-hand sides descend and ascend,
respectively, as $\alpha \in [0,\delta]$ grows, the condition is
equivalent to
$$
\frac{D}{L} > \frac{\sin2\delta}{\cos \delta-\sin\delta } = 2 \sin \delta + \frac{1}{\cos \delta}
-\frac{\cos \delta-\sin\delta(1+\sin 2\delta)}{\cos \delta(\cos
\delta-\sin\delta)}
$$
The numerator of the second ratio decays on the interval $0<\delta <
26.26^\circ$; its value at the right end $>0.1>0$ by and elementary
estimation. So this ratio is positive. Thus the inequality holds by
\eqref{cond1} and $O_3$ is disjoint with the $\delta$-extended
$\delta$-hat of $O$. For $O_4$, the conclusion is the same.
\end{document}